\documentclass[ijoc]{informs3} 

\OneAndAHalfSpacedXII 

\usepackage{natbib}
 \bibpunct[, ]{(}{)}{,}{a}{}{,}%
 \def\bibfont{\small}%
\newcommand{\secondrev}[1]{#1}

\newcommand{\thirdrev}[1]{\textcolor{black}{#1}}
\newcommand{\fourthrev}[1]{\textcolor{black}{#1}}
\newcommand{\fifthrev}[1]{\textcolor{black}{#1}}
\usepackage{algorithm}
\usepackage{algorithmicx}
\usepackage{algpseudocode}
\usepackage{subfigure}
\usepackage{enumitem}
\usepackage{amsmath}

\usepackage[table]{xcolor}

\usepackage{comment}

\newcommand{\timeSet}{\mathcal{T}}

\newcommand{\buildingSet}{\mathcal{D}}
\newcommand{\timeRoundTimeToBuilding}[1]{\tau(#1)}
\newcommand{\timeTerminalToBuilding}[1]{\tau^{1}(#1)}
\newcommand{\timeStopAtBuilding}[1]{\tau^{2}(#1)}
\newcommand{\timetBuildingToTerminal}[1]{\tau^{3}(#1)}

\newcommand{\nbuildings}{\secondrev{|\buildingSet|}}
\newcommand{\passengerStop}[1]{s(#1)}

\newcommand{\terminal}{\mathbf{T}_0}
\newcommand{\maxTime}{t^{\mathrm{max}}}

\newcommand{\passengerSet}{\mathcal{J}}
\newcommand{\nPassengers}{n}
\newcommand{\nPassengersPerBuilding}[1]{n_{#1}}

\newcommand{\passengerTimeRequest}[1]{t^r(#1)}

\newcommand{\timewindow}{T_w}

\newcommand{\cvSet}{V}
\newcommand{\nCVs}{\secondrev{|\cvSet|}}
\newcommand{\cvCap}{v^{\mathrm{cap}}}

\newcommand{\locationSet}{\mathcal{S}}
\newcommand{\location}{s}

\newcommand{\terminalnode}{\mathbf{t}}
\newcommand{\rootnode}{\mathbf{r}}

\newcommand{\set}[1]{\left\{#1\right\}}

\newcommand{\arcroot}[1]{\psi(#1)}
\newcommand{\arcterminal}[1]{\omega(#1)}
\newcommand{\nodelayer}[1]{\ell(#1)}

\newcommand{\arclayer}[1]{\ell(#1)}
\newcommand{\arccost}[1]{\eta(#1)}
\newcommand{\arcstarttime}[1]{t^0(#1)}

\newcommand{\arcdomain}[1]{\phi(#1)}

\usepackage{soul}

\newcommand{\pathSet}{\mathcal{P}}

\newcommand{\pathcost}[1]{\eta(#1)}

\newcommand{\arcTimeIndicator}[2]{\chi(#1,#2)}
\newcommand{\pathHits}[2]{k(#1,#2)}

\newcommand{\trainTrips}{\mathcal{C}}

\newcommand{\stationSet}{\mathcal{S}}

\newcommand{\rel}{\texttt{rel}}
\newcommand{\ded}{\texttt{ded}}

\newcommand{\tripTerminalTime}[1]{\tilde{t}(#1)}
\newcommand{\trainStartTime}[2]{\tilde{t}(#1,#2)}

\newcommand{\building}{d}

\newcommand{\group}{g}
\newcommand{\groups}[1]{\mathsf{g}\left(#1\right)}

\newcommand{\tottime}[1]{\textrm{tt}_{#1}}

\newcommand{\DD}{\mathsf{D}}
\newcommand{\DDnodes}{\mathsf{N}}
\newcommand{\DDarcs}{\mathsf{A}}
\newcommand{\DDnode}{\mathsf{u}}
\newcommand{\DDlayer}{\mathsf{L}}
\newcommand{\DDstate}[1]{\mathsf{s}(#1)}

\newcommand*{\permcomb}[4][0mu]{{{}^{#3}\mkern#1#2_{#4}}}
\newcommand*{\perm}[1][-3mu]{\permcomb[#1]{P}}

\usepackage{nameref}
\makeatletter
\newenvironment{taggedsubequations}[1]
{%
	\addtocounter{equation}{-1}%
	\begin{subequations}%
		\def\@currentlabel{#1}%
		\renewcommand{\theequation}{#1.\arabic{equation}}%
	}
	{\end{subequations}}
\makeatother 

\usepackage{tikz}

\usepackage{amsmath}
\usepackage{amssymb}
\usepackage{url}

\newtheorem{lem}{Lemma}
\newtheorem{thm}{Theorem}
\newtheorem{cor}{Corollary}

\newtheorem{prp}{Proposition}

\TheoremsNumberedThrough     

\EquationsNumberedThrough    

\begin{document}
\setstcolor{red}

 \RUNAUTHOR{Raghunathan et al.} 

 \RUNTITLE{Seamless Multimodal Transportation Scheduling}

 \TITLE{Seamless Multimodal Transportation Scheduling}

\ARTICLEAUTHORS{%
		\AUTHOR{Arvind U. Raghunathan}
		\AFF{Mitsubishi Electric Research Labs, Cambridge, MA, 02139, USA \\ 
			\EMAIL{raghunathan@merl.com} \URL{}}
				\AUTHOR{David Bergman}
		\AFF{Operations and Information Management, University of Connecticut,  Storrs, Connecticut 06260, USA\\ 
			\EMAIL{david.bergman@uconn.edu} \URL{}}
		\AUTHOR{John N. Hooker}
	\AFF{Tepper School of Business, Carnegie Mellon University, Pittsburgh, PA, 15213, USA \\ 
	\EMAIL{jh38@andrew.cmu.edu} \URL{}}
\AUTHOR{Thiago Serra}
\AFF{Freeman College of Management, Bucknell University, Lewisburg, PA, 17837, USA \\ 
\EMAIL{thiago.serra@bucknell.edu} \URL{}}
	\AUTHOR{Shingo Kobori}
	\AFF{Advanced Technology R\&D center, Mitsubishi Electric Corporation , Hyogo, 661-8661, Japan \\
		\EMAIL{Kobori.Shingo@cj.MitsubishiElectric.co.jp}
}} 

\ABSTRACT{%
Ride-hailing services have expanded the role of shared mobility in 
passenger 
transportation systems, creating new markets and creative planning solutions for major urban 
centers. In this paper, we consider their use for \secondrev{the first-mile or} last-mile passenger transportation in coordination 
with a mass transit service to provide a seamless multimodal transportation experience for the 
user. 
A system that provides passengers with predictable information on travel and waiting times in their 
commutes is immensely valuable.  We envision that the passengers will inform the system 
of their desired travel and arrival windows so that the system can jointly optimize the 
schedules of passengers.  
The problem we study balances minimizing travel time and the number of trips taken by the last-mile vehicles, so that long-term planning, maintenance, and environmental impact 
are all taken into account.  We focus 
on the 
case where the last-mile service 
aggregates passengers by destination.  We show that this problem is NP-hard, and  propose a decision diagram-based branch-and-price decomposition model that can solve instances of real-world size 
(10,000 passengers \secondrev{spread over an hour}, 50 last-mile destinations, 600 last-mile vehicles) in \secondrev{computational} time ($\sim 1$ minute) 
that is orders-of-magnitude faster than other methods appearing in the literature. Our experiments also indicate that \secondrev{aggregating passengers by destination on the } 
last-mile service provides high-quality solutions to more general settings. 
}%


\HISTORY{This paper was first submitted on July 2, 2019.}

\maketitle

%


\section{Introduction}
\label{sec:introduction}

Shared mobility is gradually changing how people live and interact in urban centers~\citep{CityLogistics}. 
According to McKinsey \& Company, 
the shared mobility market for China, Europe, and the United States totalled almost \$54 billion in 2016, and it is expected to grow at least 15\% annually over the next 15 years~\citep{McKinsey}. 
There is wide interest in integrating these emerging modes of transportation with public transit systems, as is indicated by the U.S. Department of Transportation~\citep{USDOT}, with interest expressed for collaboration by other key stakeholders 
such as 
companies and public sector entities signing the Shared Mobility Principles for Livable Cities~\citep{SharedMobilityWWW}. 
\secondrev{In the future, autonomous vehicles may be combined with other modes of transportation to reduce traffic congestion and the need for car ownership.}

In this paper, 
we consider the use of shared vehicles such as a bus, van,  or taxi 
in the last-mile passenger transportation, 
a particular form of transportation on demand~\citep{ToD}.  Last-mile 
transportation is defined as a service that delivers people from a hub of mass transit service to each passenger's final destination.  
Mass transit services can comprise air, boat, bus, or train. 
The last-mile service can be facilitated by 
bike~\citep{LiuJiaChe12}, car~\citep{Sha04,ThienPhd13}, autonomous pods~\citep{SheZhaZha17}, 
or personal rapid 
transit systems. 
Although last-mile may also refer to the movement of goods in 
supply chains, home-delivery systems, and telecommunications, we will restrict our attention in this paper 
exclusively to the transportation of people.  
A last-mile service expands the access of mass transit to an area wider than that defined as ``walking distance'' 
of a transportation hub.  Interest in the design and operation of last-mile services has grown tremendously 
in the past decade.  This has been 
driven primarily by three factors~\citep{Wang17}: (i) governmental push to reduce congestion and 
air pollution; (ii) increasing aging population in cities; and (iii) providing mobility for the 
differently abled and school children.

\begin{figure}[ht]
	\begin{center}
		\includegraphics[scale=0.25]{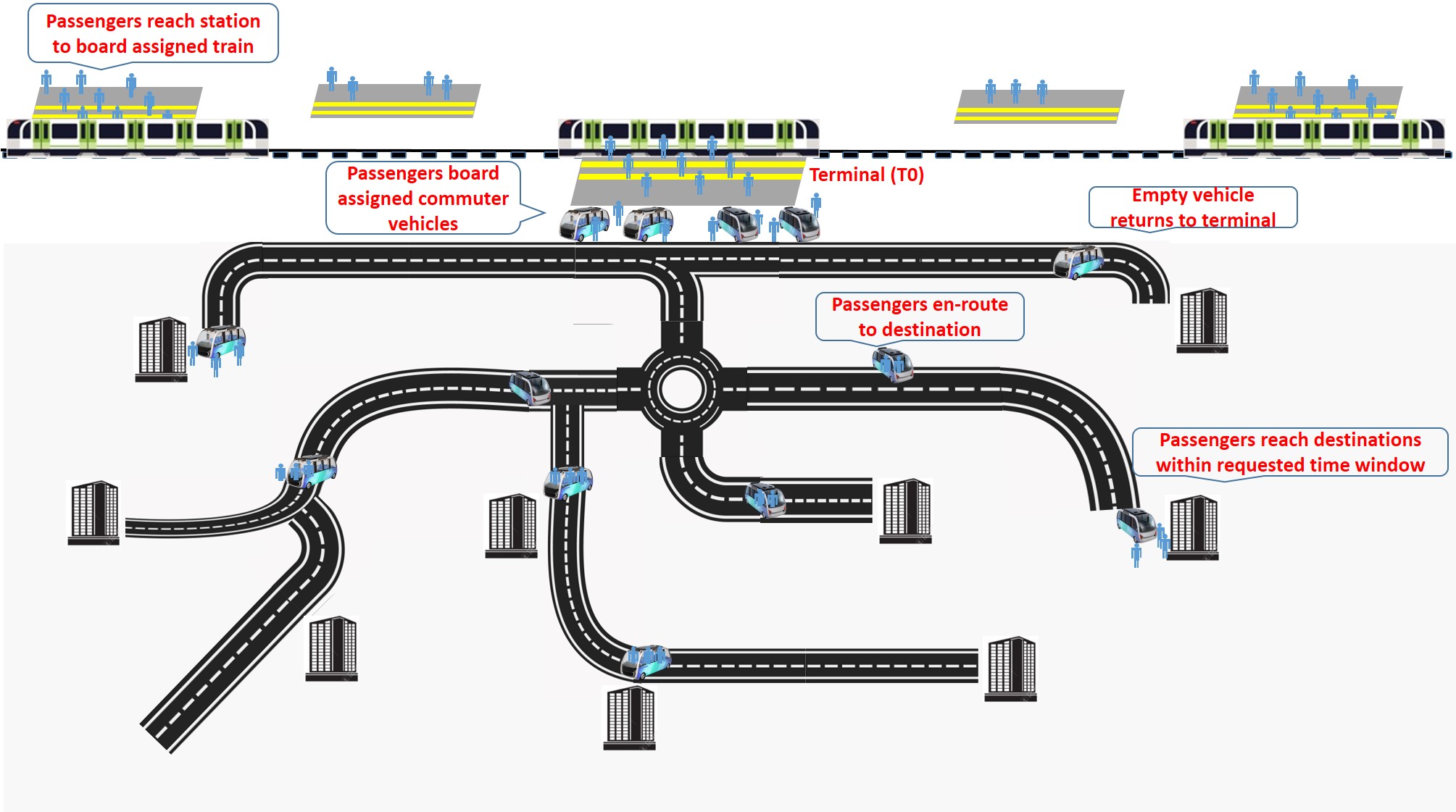}
		\caption{Schematic of an integrated last-mile system (images licensed from \textcolor{blue}{\texttt{shutterstock.com}} and \textcolor{blue}{\texttt{alamy.com}}).}
		\label{fig:ProbDes}
	\end{center}
\end{figure}
Figure~\ref{fig:ProbDes} shows a typical scenario for the operation of mass transit in conjunction 
with a last-mile service.  All passengers start their journey from mass transit stations 
served, for example, by a train, and request transportation to their destinations  
within a time-window.  Destinations can represent offices in specific buildings or 
a location, such as a bus stop, that serves multiple buildings that are easily accessed on foot.  
These destinations can be accessed 
by paths that may be exclusive to last-mile vehicles or by means of existing road infrastructure.  
For convenience, we will refer to the vehicle providing last-mile service as a 
\emph{commuter vehicle} (CV).  The CVs are typically parked at a terminal ($\terminal$) at 
which passengers arrive from mass transit services and 
proceed to their respective destination buildings by sharing a ride in a CV.  
The last-mile service may represent the morning commute to the office or the evening commute back to 
residences.   Once all passengers are delivered to their \fifthrev{common destination}, \secondrev{it is assumed that} the CVs return back to 
the terminal for subsequent trips. \fifthrev{Namely, we make the following assumptions:
(1) all passengers of a trip share the same destination; and
(2) CVs return to their depot after a trip.}

\secondrev{The particular system that we describe is currently being planned 
for the city of Milton-Keynes in the United Kingdom.  The project called UK 
Autodrive (\url{http://www.ukautodrive.com/}) trialled the use of autonomous pods for the 
passenger transportation within the city.  In particular, the city officials in Milton-Keynes are 
looking to deploy this system for transport of passengers between the Milton-Keynes train station and 
the city limits\footnote{\url{http://www.ukautodrive.com/pods-provide-a-first-last-mile-solution-in-milton-keynes/}}.}

Such an application is timely and well suited for integrated prescriptive analytics. 
This type of system is already realizable in practice by integrating with ride-hailing services such 
as conventional taxi services, Uber, and Lyft.  
Moreover, the peak  use of the system is expected to coincide 
with morning and afternoon work commutes, which can be made predictable by design: 
work commuters often know in advance when they would like to arrive at work 
and would be willing to provide such information in advance for better service.  \secondrev{Note that this assumption on operation mode does not preclude a real-time decision making process.  One could adopt a rolling horizon approach, fixing decisions for previously scheduled passengers and re-optimizing for newly arrived requests as often as needed.}
Finally, note that the \emph{first-mile} operation, wherein the passengers first ride on CVs to reach 
a hub of a mass transit service, can easily be accommodated in an analogous manner\secondrev{, as long as there are shared pickup locations that will serve as hubs}.

\begin{figure}[h]
	\begin{center}
		\includegraphics[scale=0.3]{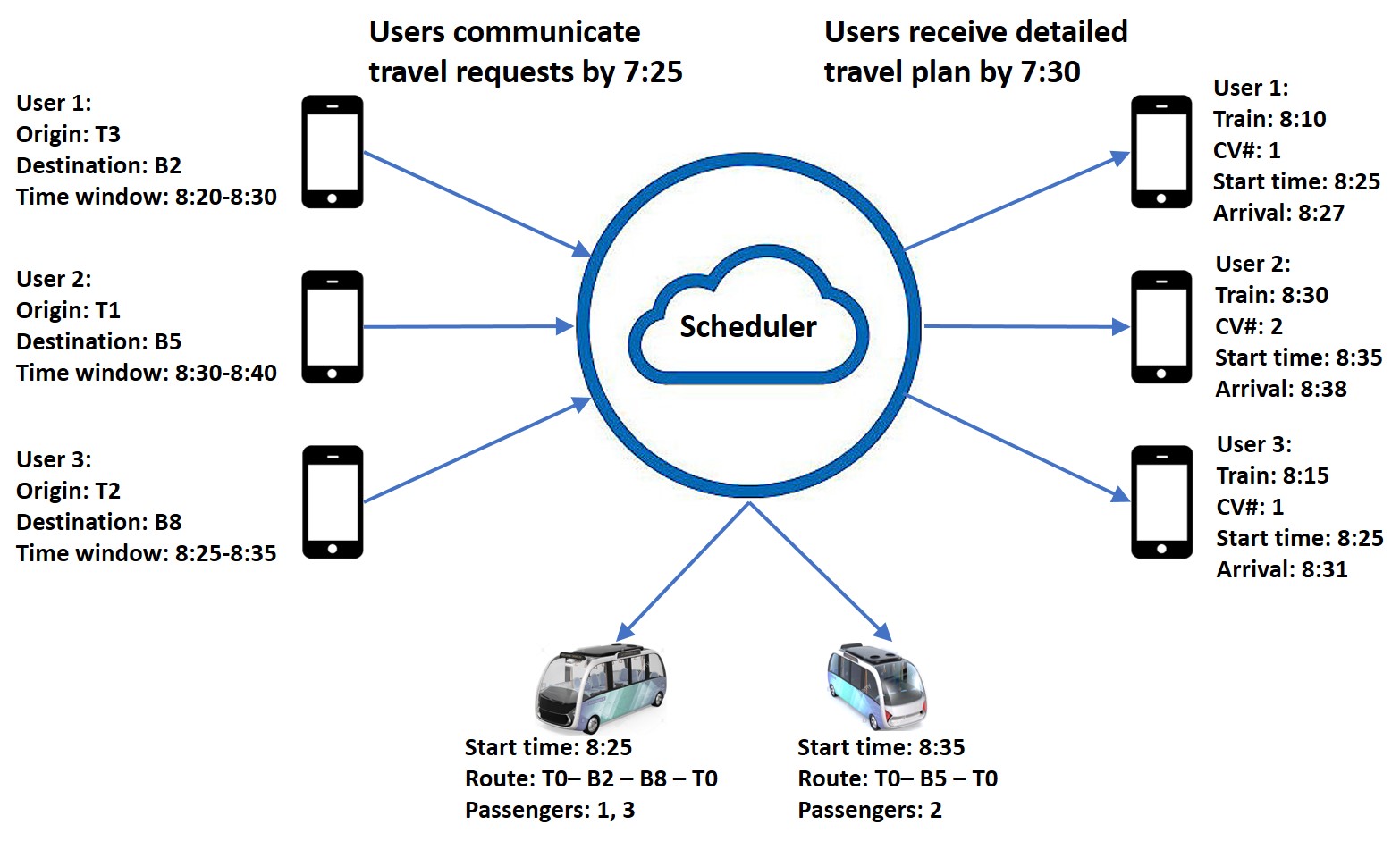}
		\caption{Schematic representation of the interaction between the passengers in the system, the scheduler, and the CVs (images licensed from \textcolor{blue}{\texttt{shutterstock.com}} 
		 and \textcolor{blue}{\texttt{alamy.com}}) }
		\label{fig:UserInteracn}
	\end{center}
\end{figure}
Therefore, we envision an operational scenario such as in Figure~\ref{fig:UserInteracn}, where passengers 
indicate their station of origin, destination, and the desired time-window of arrival at 
the destination.  
\secondrev{This information is assumed to be available to the scheduler only slightly in advance of 
scheduling decisions.}  For instance, consider the situation where the destinations represent buildings 
with offices and the passengers  enter requests through a smartphone app.
\secondrev{Once all requests have been received for a certain time-period,} the scheduler 
determines for each user: (i) the mass transit trip to board at their station of origin; (ii) the CV to 
board at $\terminal$; (iii) the time 
the CV will depart from $\terminal$; and (iv) the time of arrival at the destination.  The scheduler also 
communicates to the different CVs the routes, start times, and list of passengers.   The 
choice of route determines the times of arrival of passengers at their destination as well as total travel time.

This application characterizes what we denote as the Integrated Last-Mile Transportation Problem 
(ILMTP).   The ILMTP is defined as the problem of scheduling passengers jointly on mass transit and 
last-mile services so that the passengers reach their destinations within specified time-windows. 
We propose to minimize a linear combination of the total transit time for all passengers and 
the number of trips required.  The former quantity captures the quality of service provided; the 
latter addresses fuel consumption, long-run operational costs, and environmental considerations.    
Transit time includes the time spent traveling in both transportation modes and waiting between both services.  
The solution of the ILMTP also specifies the set of passengers that share a ride in a CV.  The time 
spent by the passengers in the CV also depends on the co-passengers\secondrev{, where we focus on a variant simplified as discussed below}.  

\subsection{Contributions of this paper}

\thirdrev{This paper proposes algorithms that are scalable to multimodal 
transportation systems consisting of a central depot connecting a mass transit system with multiple small-capacity vehicles. Our experimental evaluation on instances generated to emulate what might be realized in real-world settings indicates 
that the algorithms developed can obtain 
optimal solutions in one minute to problems with 
10,000 passengers, 50 destinations, and 600 vehicles that commonly utilized techniques are 
not able to solve in reasonable time.  
We also obtain high-quality solutions to more general problems than the one we focus on, 
finding heuristic solutions in seconds to problem instances for which other approaches in the 
literature cannot find a feasible solution in three hours.  }


More specifically, we extend results first described in~\cite{iflmpicaps18}, where a 
single-destination-per-trip (SDPT) assumption is also imposed. 
We note that this assumption can be regarded as a technical constraint of the business model. 
First, users do not want to be delayed due to other passengers leaving the car in previous stops, 
especially if the distance between these stops is walkable. 
Second, the literature on last-mile transportation regards the use of existing destinations, such as 
bus stops, as an effective aggregator of individual destinations~\citep{MeetingPoints,maheo18}.  
We can observe some instances of this logic in practice, such as the new Uber Express POOL service~\citep{UberX}, which embodies 
these considerations by creating pooled pickup and drop-off locations.

{Aligned with \citep{iflmpicaps18}, we make the following three assumptions on the operation of the system: 
\begin{itemize}
\item The trains serves station along the line and have a deterministic travel between the stations.  
\item Each CV trip is a round trip from the terminal $\terminal$ to one of the buildings in $\buildingSet$. 
\item The passengers have identical time window on arrival at destination. 
\end{itemize}
The above assumptions imply that the ILMTP is modeled and solved for each terminal station independently when applied to a realistic transportation network.}  

The main contributions of the paper can be summarized as follows:
\begin{itemize}
	\item \textbf{Computational Complexity:} We show that the ILMTP with the SDPT assumption, 
	 which, for simplicity, we refer to simply as the ILMTP-SD unless otherwise noted, is 
	\emph{NP-Hard}.  This shows that the simplifying assumption does not render the problem 
	computationally easier to solve and developing an efficient algorithm is necessary.
	\item \textbf{Structure of Optimal Solutions:} The optimal solutions to the ILMTP-SD are 
	shown to satisfy an ordering property.
	This result is a restatement of~\cite{iflmpicaps18}[Theorem~1], 
	 included for completeness. 
We also show that the number of solutions to the ILMTP-SD grows exponentially with the number of passengers.
	\item \textbf{Decision Diagram (DD)-based Algorithm:} Based on the structure of 
	optimal solutions, we describe a novel optimization 
	algorithm based on a DD representation of the space of solutions to the ILMTP-SD.  
	Our approach builds on a state-space 	
	decompositions \citep{BerCir16,BerCir17}, which is an outgrowth of the growing body of literature 
	on DD-based optimization \citep{BerHoeHoo11,bergman2016decision}.  
	\item \textbf{Branch-and-price DD decomposition:} We improve the performance of the proposed DD algorithm 
         by developing a branch-and-price scheme~\citep{BnP} through which columns are generated from paths in the DDs. 
         This scheme could be easily adapted to other problems in which such a DD-based algorithm is devised.
	\item \textbf{Numerical Evaluation}:  An experimental evaluation indicates that the 
	proposed model is orders-of-magnitude faster than existing techniques. 
	We also report the performance of the proposed approach on instances where the 
	assumptions on the ordering property are violated, namely (i) mass transit service also includes 
	express trains and (ii) CVs 
can make multiple stops in the last-mile.  
	We show that even in these settings the potential loss of optimality incurred by the SDPT 
	assumption is insignificant for all practical purposes, in that existing models are often unable 
	to find a single feasible solution in hours while the algorithm we propose prove optimality 
	for the restricted version in seconds.
	\secondrev{Thus, existing approaches are inappropriate for real-time use in the ILMTP, especially when requests are received without much advance notice. The proposed method is capable of optimizing this problem and is appropriate as a real-time heuristic.}
\end{itemize}




\subsection{Related work}

The ILMTP can be broadly viewed as an instance of routing and scheduling with 
time-windows.  We survey the relevant literature and  describe the key differences with this paper.

The literature on last-mile transportation is mostly focused on last-mile service, 
without much consideration to the mass transit system.  Seminal work in this area dates back to the 1960s, focused 
mostly on freight transportation (see~\citeauthor{Wang17}~(\citeyear{Wang17}) for a discussion).  
\cite{Wang17} is the first work to consider routing and scheduling in the last-mile, where the 
minimization of total travel time is considered and the author proposes a heuristic 
approach.  The ILMTP is a strict generalization
of~\cite{Wang17}, in that we consider time-windows for arrival and scheduling on the mass transit 
service. 
An alternative approach to exact last-mile integration consists of simulation-based approaches~\citep{Simulations1,Simulations2} as well as last-mile dispatching policies~\citep{Dispatching1,Lau2018}. 
This is complemented by another stream of work on uncertainty of mass transit service~\citep{TransitUncertainty1,TransitUncertainty2}. 
More recently,~\cite{maheo18} approached the design of a public transit system 
that includes multiple modes of transportation, however the authors did not consider 
scheduling aspects.  
This paper expands this literature with an optimal approach to scheduling 
passengers in the last-mile\st{.}\secondrev{,} 
and complements the work of  \cite{maheo18}  by focusing exclusively on scheduling. 

Personal Rapid Transit (PRT) has
similarities to the last-mile problem.
Research has been conducted on PRT system control frameworks~\citep{And98}, financial
assessments~\citep{BlyTey05,Beretal11}, performance approximations~\citep{LeeHamDav09,LeeHamWil10}, 
and case studies~\citep{MulSgo11}.  However, none of these papers have addressed last-mile 
operational issues.

On the other hand, a large body of research has been devoted to Demand Responsive Transit (DRT), which is another type of
on-demand service. Some papers focus on DRT concept discussions, practical implementation, and 
assessment of simulations in case 
studies~\citep{BraNelWri04,Horn02,MagNel03,PalDesAbd04,QuaDesOrd08}.
Models have been developed to assist in system design and regulation
(e.g.,~\cite{Dag78},~\cite{DiaDesXia06},~\cite{WilHen80}).  Routing options in specific contexts have 
also been considered~\citep{Cheetal12,Horn02b}. The ILMTP can be viewed as a specific variant of a 
broadly defined DRT concept\textemdash  namely, a DRT system that addresses 
last-mile service requests with batch passenger demand and a shared passenger origin.  
The same applies to ride-sharing models more generally~\citep{DynamicRideSharing}. 
Unlike 
most papers in the DRT literature, we focus on scheduling the mass transit service and 
the last-mile optimization from an operational perspective.

A much broader stream of related work consists of 
vehicle routing problems (VRPs), 
which have long been studied and comprise a large body of literature.  
The VRP with time windows (VRPTW) has been the subject of intensive study, 
with many heuristic and exact optimization approaches suggested in the literature~\citep{vrpbook}. The dial-a-ride problem (DARP) 
and related variations have also 
been extensively investigated~\citep{CorLap07,Jawetal86,LeiLapGuo12}.  
As argued by~\cite{Wang17}, the VRPTW focuses on reducing operating costs while the ILMTP 
aims to improve the level-of-service by minimizing total passenger transit time, 
especially with single-destination routes. 
The typical size of the problems that can be solved to optimality for the VRPTW and DARP 
are on the order of 
hundreds 
of requests and 
dozens 
of vehicles. This is far smaller than the 
size of the instances that we solve to optimality in this paper.  Although solving the VRPTW and DARP optimally
is difficult for large-size instances, good heuristics exist for 
these problems, 
which include
the Savings algorithm 
\citep{ClarkeWright64}, 
its variants, and insertion heuristics~\citep{Vigo96,SalhiNagy99,CampSavel04}.  
In a recent work, \cite{RileyLegrainHentenryck19} propose a column generation method for scheduling 
passengers on ride-sharing vehicles which scales to 30,000 requests per hour, but does not consider the aspect of mass transit service and 
time-windows.

Finally, 
the computational approach introduced in this paper calls for modeling problems using disjoint and connected 
DDs. 
DD-based optimization is an emerging field within computational optimization \citep{AndHadHooTie07,BerHoeHoo11,Gange2011,CirHoe13,bergman2016decision,DBLP:conf/aaai/PerezR17}.  The idea used in this paper is to model a problem with a set of 
DDs, from which a solution is obtained by a network flow reformulation based on such DDs~\citep{BerCir16,BerCir17}.  This paper proposes a path-based model for solving the resulting network-flow reformulation, 
which is similar to the approach 
by \cite{MorSew13} for the vertex coloring problem. 
To the best of our knowledge, this paper is the 
first application to multi-valued DDs.

In summary, the ILMTP has the following features that distinguishes it from previous studies 
in the literature: 
\begin{itemize}
\item joint scheduling of passengers on mass transit systems and last-mile services; 
\item consideration of time-windows on arrival at destination; 
\item common last-mile origin (which is also the 
vehicle depot); and  
\item weighted minimization of the total passenger transit time and number of CV trips.
\end{itemize}

The ILMTP models real-world transportation systems that are prevalent 
across the globe, and this paper provides a mechanism for optimal operational decisions.

\section{Problem Description and Mathematical Formulation}\label{sec:problemDescription}

%
%

In this section we provide a formulation of the ILMTP-SD.  Prior to that, we describe the  
different elements in the ILMTP-SD such as the mass transit system, last-mile vehicles, destinations, 
passenger requests and associated parameters such as the travel time associated with the 
transportation services, and time windows for arrival.

\noindent\textbf{Mass transit system:} For the sake of convenience, we will assume that the mass 
transit is a train system.  Let $\terminal$ be the \emph{terminal station} that links a mass transit 
system with a last-mile service system.  
The mass transit system is described by a set of \emph{trips}, denoted by $\trainTrips$. 
Each trip originates at a station in  set $\locationSet$ and ends at $\terminal$.  The trips are regular in the sense that the train stops at all stations in $\locationSet$ sequentially, with $\terminal$ as the last stop of each trip. The time a trip $c$ leaves station $\location \in \locationSet$ is 
$\trainStartTime{c}{\location}$ and the time it arrives to the terminal is $\tripTerminalTime{c}$. This paper is only concerned with the moving  portion of the mass transit commute and so it assumes that each passenger arrives at the station of origin  at the time that the mass transit trip is departing that location.  The time that a passenger waits in such stations is not of concern in our objective or constraints. 

\noindent\textbf{Destinations:} Let $\buildingSet$ be the set of destinations where the CVs make stops,
where we assume $\terminal \in \buildingSet$.    
For each destination $d \in \buildingSet$,  let $\timeRoundTimeToBuilding{d}$ be the total time it takes a CV to depart $\terminal$, travel to $d$ (denoted by $\timeTerminalToBuilding{d}$), stop at $d$ for passengers to disembark (denoted by $\timeStopAtBuilding{d}$), and return to $\terminal$ (denoted by $\timetBuildingToTerminal{d}$). Therefore, $\timeRoundTimeToBuilding{d} = \timeTerminalToBuilding{d} + \timeStopAtBuilding{d} + \timetBuildingToTerminal{d}$. Let $\timeSet := \set{1, \ldots, \maxTime}$ be an index set of the operation times of both systems.  We assume that the time required to board 
passengers into the CVs is incorporated in $\timeTerminalToBuilding{d}$.  
For simplicity, the boarding time is independent of the number of passengers.  A
passenger arrives to a destination $\timeTerminalToBuilding{d}$ time units 
after departing from the terminal.

\noindent\textbf{Last-mile system:} Let $\cvSet$ be the set of CVs.
Let $\cvCap$ denote the \secondrev{maximum} number of passengers that can be assigned to a single CV trip.   Each CV trip consists of a set of passengers boarding the CV, traveling from $\terminal$ to a destination $d \in \buildingSet$, and then returning back to $\terminal$.  
Passengers on a common CV trip must request transportation to a common building. We also assume that each CV must be back at the terminal by time $\maxTime$.

\noindent\textbf{Passengers:} Let $\passengerSet$ be the set of passengers.  Each passenger $j  \in \passengerSet$ requests 
\secondrev{train} from a station $\passengerStop{j} \in \locationSet$ to $\terminal$, and then by CV \secondrev{from $\terminal$} to destination $d(j) \in \buildingSet$, to arrive at time $\passengerTimeRequest{j}$.  The set of passengers that request service to destination $d$ is denoted by $\passengerSet(d)$.  Let $\nPassengers := |\passengerSet|$ and $\nPassengersPerBuilding{d} := \left| \passengerSet(d) \right|$. Each passenger $j \in \passengerSet$ must arrive to $d(j)$ between $\passengerTimeRequest{j} - \timewindow$ and  $\passengerTimeRequest{j} + \timewindow$.

\noindent\textbf{Problem Statement:} 
The ILMTP-SD is the problem of assigning train trips and CVs to each passenger so that the total transit time and the number of CV trips utilized is minimized.   A solution therefore consists of a partition $\mathsf{g} = \set{g_1, \ldots, g_\gamma}$ of $\passengerSet$, with each group $g_l$ associated with a departure time $t^{\mathsf{g}}_l$, for $l = 1, \ldots, \gamma$, which indicates the time the CV carrying the passengers in $g_l$ departs $\terminal$, satisfying all request time and operational constraints.  For any passenger $j \in \passengerSet$, let $\mathsf{g}(j)$ be the group in $\mathsf{g}$ that $j$ belongs to.

To balance the potentially conflicting objectives, the objective function we consider is a convex combination of objective terms, defined by $\alpha$, with $0 \leq \alpha \leq 1$. Hence, we balance these objectives by using $\alpha$ times the waiting time plus $(1-\alpha)$ times the number of CVs, 
which is therefore used as our objective function, represented as $f(\alpha)$.  

\subsection{IP Model}

In this section we present an improved IP model for the ILMTP-SD, which is based on the one from~\cite{iflmpicaps18}.  
\secondrev{Let 
$\timeSet_j : = [\passengerTimeRequest{j} - \timewindow - \timeTerminalToBuilding{\building(j)},\passengerTimeRequest{j} + \timewindow-\timeTerminalToBuilding{\building(j)}]$ 
represent the possible times at which a 
passenger can depart on a CV from $\terminal$.  
The total travel time for each passenger $j$ that travels on a 
CV that departs from $\terminal$ at time $t \in \timeSet_j$ is 
	\begin{equation*}
	\tottime{j,t} := \left( t+\timeTerminalToBuilding{\building(j)} - \max\limits_{c \in \trainTrips : \tripTerminalTime{c} \leq t} \trainStartTime{c}{\passengerStop{j}} \right).
	\end{equation*}
which is composed of travel time on the train, waiting time 
at $\terminal$, and CV travel time. Note that the passenger is assigned to the latest train that allows to reach $\terminal$ prior to departure on the CV. }
The variables in our model are as follows:

\begin{itemize}


	\item $z_{j,t}$: indicator if passenger $j \in \passengerSet$ leaves $\terminal$ at time $t \in \timeSet_j$
	\item $n_t$: number of CVs parked in $\terminal$ at time $t \in \timeSet$
	\item $n_{d,t}$: number of CVs assigned to destination $d \in \buildingSet$ to depart $\terminal$ at time $t \in \timeSet$
\end{itemize}
An optimization model for ILMTP-SD is as follows:
\secondrev{
	\begin{taggedsubequations}{IP}
		\label{eqn:IPmodel}
		\begin{align}
		 &	\min && f(\alpha) = \alpha \cdot \sum_{j \in \passengerSet} \sum\limits_{t \in \timeSet_j}\tottime{j,t} \cdot z_{j,t} + (1-\alpha) \cdot  \sum_{d \in \buildingSet} \sum_{t \in \timeSet}  n_{d,t}
		 \label{ip1} 
		 \\
		 & \textnormal{s.t. }
		 &&
		 \sum_{t \in \timeSet_j} z_{j,t} = 1, 
		 &&
		 \forall j \in \passengerSet 
		 \label{ip3}
		 \\ 
		 &&&
		 	n_{t}
			=
			n_{t-1}
			+ 
			\sum_{d \in \buildingSet}
			n_{d,t-\timeRoundTimeToBuilding{d}}
			-
			\sum_{d \in \buildingSet}
			n_{d,t}
			&&
			\forall t \in \timeSet
			\label{ip6}
		 \\
		 &&&
				 \cvCap \cdot 
		 \left(
		  n_{d,t} - 1
		  \right) + 1 \leq  \sum_{j \in \passengerSet(d)}
		 z_{j,t} 
		 \leq 
			\cvCap \cdot n_{d,t}
			&&
			\forall d \in \buildingSet, 
			\forall t \in \timeSet
			\label{ip7}
		 \\
		 &&&
		 z_{j,t} \in \set{0,1},
		 &&
		 \forall j \in \passengerSet, \forall t \in \timeSet 
		 \label{ip11} \\
		 &&&
		 n_t \geq 0,
		 && 
		 \forall t \in \timeSet
		 \label{ip12}
		 \\
		 &&&
		 n_{d,t} \geq 0, && \forall d \in \buildingSet,  
		 \forall t \in \timeSet
		 \label{ip13} \\
		 &&&
		 n_0 = \nCVs. 
		 \label{ip14}
		\end{align}
 	\end{taggedsubequations}
}

The objective function, parametrized by $\alpha \in [0,1]$, 
balances the sum of the total travel time of all passengers with the number of CV trips that take place over the planning horizon.   
This objective function generalizes the one in~\cite{iflmpicaps18} by also including the number of CV trips as an element of the objective function.  The smaller $\alpha$ is, the more emphasis is placed on minimizing the number of trips a CV takes. Fewer CV  trips results in fewer maintenance tasks as well as lower emissions, which is critical for long-term planning and for minimizing environmental impact. In the case of the ILMTP-SD 
with regular train trips, 
only the waiting time at the terminal varies across solutions taking the shortest route to each destination. 
Nevertheless, we use total travel time for consistency in the section on experiments where we relax SDPT and the regular train times assumption. 

The 
\secondrev{first} set of constraints in
(\ref{ip3}) ensure each passenger is assigned to one
one CV trip.  
Constraints~(\ref{ip6}) through~(\ref{ip7}) 
bookmark the number of CVs in use at any given time.   
In contrast to previous formulations, 
the constant value 1 on constraint (\ref{ip7}) prevents empty CVs in feasible solutions. 
Constraints~(\ref{ip11}) through~(\ref{ip13}) enforce bounds, binary restrictions, and initial conditions, as necessary. Note that~\eqref{eqn:IPmodel} can be extended to handle time-dependent 
travel times on the CV by replacing the occurrence of 
$\timeTerminalToBuilding{d}$, $\timeRoundTimeToBuilding{d}$ in~\eqref{eqn:IPmodel}
with CV travel times that are dependent on the time of departure from $\terminal$.

\section{Complexity}
\label{sec:complexity}

%

It is known that generalizations of the ILMTP are NP-hard~\citep{iflmpicaps18}. 
We show in this section that the ILMTP-SD is at least as hard, 
and in fact a much simpler version of the ILMTP-SD with a single mass transit service and 
CVs of unitary capacity is sufficient to define an NP-hard problem, with proof in Appendix Section~\ref{sec:proofCompleixty}.

\begin{thm}
	\label{thm:complexity}
	Deciding feasibility of an instances of the ILMTP-SD is NP-complete.  
\end{thm}


\begin{cor}
ILMTP-SD is NP-hard.
\end{cor}


\section{Structure of Optimal Solution to ILMTP-SD}
\label{sec:solnStructure}

In this section we restate a result from~\cite{iflmpicaps18} uncovering a structural property of optimal solutions to the ILMTP-SD. \secondrev{This structural result enables us to construct a recursive model for solutions to the ILMTP-SD as will be discussed in \S~\ref{sec:singleDestBDD}.}
We also prove a exponential lower bound  on the number of 
solutions to the ILMTP-SD in Theorem~\ref{thm:counting}.  
Note that the assumptions stated in~\cite{iflmpicaps18} hold in the 
present context (see \S~\ref{sec:problemDescription}).  
\begin{thm}\label{prop:optSolutionStructure}
For all $d \in \buildingSet$, let $\passengerSet(d) = \set{j^d_1, \ldots, j^d_{\nPassengersPerBuilding{d}}}$ represent a partitioning of $\passengerSet$ by destination and let $\passengerSet(d)$ 
be ordered so that, for $1 \leq i \leq \nPassengersPerBuilding{d}-1$, $\passengerTimeRequest{j^d_i} \leq \passengerTimeRequest{j^d_{i+1}}$.  There exists an optimal solution for which there are no triples of passengers $j^d_{i_1}, j^d_{i_2}, j^d_{i_3}$ with $i_1 < i_2 < i_3$ for which $j^d_{i_1}$ and $j^d_{i_3}$ share a common CV trip without $j^d_{i_2}$. 
\end{thm}
Theorem~\ref{prop:optSolutionStructure} indicates that one need only search for solutions which group passengers in CV trips sequentially by order of requested arrival times.   We will exploit this result in order to create compact DDs for each destination.  The proof of  Theorem~\ref{prop:optSolutionStructure} follows from the proof of Theorem 1 from~\cite{iflmpicaps18}, where we note that in the exchange argument, no CV trips are added or deleted. We provide the proof in 
Appendix~\ref{app:theorem2}.

\subsection{Exponential Lower Bound on Number of Solutions}

By Theorem~\ref{prop:optSolutionStructure}, 
the search for optimal solution\secondrev{s} can be restricted to groups of passengers that are consecutive when ordered by deadlines.  Theorem~\ref{thm:complexity} shows that the problem remains NP-hard even over these solutions. Theorem~\ref{thm:counting} relates that number of partitions of passengers going to the same destination to the Fibonacci series establishing a exponential lower bound on the number of feasible solutions even in this restricted set.

\begin{thm}\label{thm:counting}
	Let $\phi(n)$ be the number of partitions of $n$ passengers into groupings, each containing passengers with consecutive deadlines and going to a common destination.  If the time windows and requested arrivals times are such that every pair of consecutive passengers can travel with one another on a common CV, then the number of partitions of passengers into groups is bounded below by the $(n+1)$st Fibonacci number $\mathsf{F}(n)$, and is hence exponential in $n$ ($\phi(n) \sim O(1.6^n)$).
\end{thm}

\proof{Proof.}
Consider the case when $\cvCap = 2$. For $n \geq 3$, we have \st{htat}\secondrev{that}
$ \phi(n) = \phi(n-1) + \phi(n-2).$
This can be shown by conditioning on whether or not the last passenger travels alone or with the penultimate passenger.  In the first case, of  traveling alone,  the number of partitions of the remaining set of passengers is $\phi(n-1)$.  In the second case, of traveling with 
another passenger, the number of partitions of the other passengers is  $\phi(n-2)$.  Since $\phi(1) = 1$ and $\phi(2) = 2$, the result follows. 
For larger $\cvCap$, the recursion is written
\[ \phi(n) = \phi(n-1) + \phi(n-2) + \cdots + \phi(n - \cvCap), \]
which is bounded from below by $\phi(n-1) + \phi(n-2)$.
\Halmos
\endproof


\section{State-Space Decomposition}
\label{sec:ssDecomp}

In this section we discuss the modeling of the ILMTP-SD through \emph{decision diagram (DD) decomposition} \citep{BerCir16,BerCir17}, which relates to \emph{state-space decompositions} using dynamic programming \citep{Bertsekas1999,Bertsekas2012}. In particular, we show how one can model every possible single-destination CV trip through a compact DD.  We then describe how the DDs can be concurrently optimized over through a network-flow reformulation with channeling constraints, which provides a novel and computationally advantageous remodeling of the ILMTP-SD relying on the structural property in Theorem~\ref{prop:optSolutionStructure}.
Section~\ref{sec:BDDgeneration} describes how such a collection of DDs can be efficiently constructed,
and Section~\ref{sec:NFR} shows how to jointly use the constructed DDs to find an optimal solution for the problem.

\subsection{Single destination BDD}
\label{sec:BDDgeneration}
\label{sec:singleDestBDD}

For each destination $d \in \buildingSet$ we construct a 
DD that encodes every possible partition of $\passengerSet(d)$ into 
CV trips through paths.  A DD is an acyclic digraph 
with nodes that are split among layers and arcs that connect 
nodes in different layers. DDs have close connection to the 
stage-state representation in dynamic 
programming~\citep{Bertsekas2012}. 
The layers in the DD can be identified with stages, the nodes in 
each layer represent different states and the arcs in the DD 
represent feasible transitions between states in successive stages. 

Theorem~1 shows that it is 
sufficient to only consider groups of passengers with consecutive 
deadlines. Let the passengers in $\passengerSet(d)$ be ordered in 
nondecreasing order of $\passengerTimeRequest{j}$ as 
$j^d_1, \ldots, j^d_{\nPassengersPerBuilding{d}}$. 
In order to encode the 
set of all partitions of $\passengerSet(d)$ into 
CV trips, the key decision for a pair of consecutive 
passengers $j^d_i, j^d_{i+1}$ is whether they ride together in a 
CV trip. The decision of $j^d_i$ sharing a ride in a CV with 
$j^d_{i+1}$ is predicated on (a) prior decisions (other passengers 
$j^d_{i-1}, j^d_{i-2},...,j^d_{i-\cvCap+1}$ already sharing a 
vehicle with $j^d_i$) since this affects the availability of a seat 
for $j^d_{i+1}$ and (b) the time windows of the passengers should have 
non-empty intersection indicating that there exists a common time 
of arrival. 

Accordingly, the layer in the DD (stage in DP) will denote the index 
of the passengers whose decisions have already been considered, the  
nodes in the DD (states in DP) will denote the number of occupants 
in the CV, and finally the 
arcs between the layers in the DD will represent the decision  
of a passenger represented by the layer to ride 
together with the next passenger in $\passengerSet(d)$.  
The addition of the arcs will satisfy the stated conditions in 
(a) and (b).  Further, we associate additional information on 
start time and the objective contribution from CV trip on the arcs. 
Additionally, each path establishes the departure time of each CV 
and the total contribution to the objective function of the 
passengers in $\passengerSet(\building)$ 
given the partition prescribed by the path.  We make these 
precise in the following.

Formally, we construct, for every \secondrev{destination }$d \in \buildingSet$, a \st{DD} \secondrev{diagram} $\DD^d$ , which is a layer-acyclic digraph 
$\DD^d = (\DDnodes^d,\DDarcs^d)$.  $\DDnodes^d$ is partitioned into $\nPassengersPerBuilding{d} + 1$ ordered layers $\DDlayer^d_1, \DDlayer^d_2, \ldots, \DDlayer^d_{\nPassengersPerBuilding{d}+1}$ where $\nPassengersPerBuilding{d} = 
|\passengerSet(d)|$.  Layer $\DDlayer^d_1 = \set{\rootnode^d}$ and layer $\DDlayer^d_{\nPassengersPerBuilding{d} + 1} = \set{\terminalnode^d}$ consist of one node each; the \emph{root} and \emph{terminal}, respectively. The layers $i = 1,...,\nPassengersPerBuilding{d}$ correspond to the passenger $j^d_i$. 
The \emph{layer} of node $\DDnode \in \DDlayer^d_i$ is defined as 
$\nodelayer{\DDnode} = i$. Each arc $a \in \DDarcs^d$ is directed from 
its \emph{arc-root} $\arcroot{a}$ (a node in layer $i$) to its 
\emph{arc-terminal} $\arcterminal{a}$ (a node in layer $i+1$), 
\emph{i.e.} 
$\nodelayer{\arcroot{a}} = \nodelayer{\arcterminal{a}} - 1$.  
We denote the \emph{arc-layer} of $a$ as 
$\arclayer{a} := \nodelayer{\arcroot{a}}$.  It is assumed 
that every maximal arc-directed path connects $\rootnode^d$ to 
$\terminalnode^d$.
Each node $\DDnode$ is associated with a \emph{state} 
$\DDstate{\DDnode}$ that defines the passengers aboard a CV trip.  
Since the capacity of the CV is $\cvCap$ the possible 
states are $0,\ldots,(\cvCap-1)$. An arc between the 
$\DDnode \in \DDlayer^d_i$ and $\DDnode' \in \DDlayer^d_{i+1}$ 
exists if (i) $\DDstate{\DDnode'} = \DDstate{\DDnode}+1$ or 
$\DDstate{\DDnode'} = 0$ and 
(ii) if $\DDstate{\DDnode'} \neq 0$, 
$\cap_{k \in \{0,\ldots,\DDstate{\DDnode'}\}} 
[\passengerTimeRequest{j^d_{i+1-k}} - \timewindow,\passengerTimeRequest{j^d_{i+1-k}} + 
\timewindow] \neq \emptyset$.
It is easy to verify that (i) and (ii) imply the satisfaction of 
(a) and (b) stated earlier in the section. 
We distinguish between the arcs that connect to the node 
with state $0$ in the next layer and the remaining arcs.  
The former set of arcs are called \emph{one-arcs} and the 
rest of the arcs are called the \emph{zero-arcs}, 
indicated by 
\secondrev{one-arcs, $\arcdomain{a} = 1$ and zero-arcs, $\arcdomain{a} = 0$}. 
In fact, we duplicate the one-arcs as many times as there are 
time instances in the set 
$\cap_{k \in \{0,\ldots,\DDstate{\DDnode}\}} 
[\passengerTimeRequest{j^d_{i-k}} - \timewindow,\passengerTimeRequest{j^d_{i-k}} + 
\timewindow]$.

Every one-arc $a$ corresponds to a group
 $\group(a) = \left\{ 
	j^d_{\arclayer{a} - \DDstate{\arcroot{a}}} , j^d_{\arclayer{a} - \DDstate{\arcroot{a}} + 1}, \ldots , j^d_{\arclayer{a}}\right\}$, i.e, the set of contiguously indexed $\DDstate{\arcroot{a}}+1$ 
	passengers ending in index $\arclayer{a}$.
A one-arc stores an \emph{arc-cost} $\arccost{a}$ and an \emph{arc-start-time} $\arcstarttime{a}$.  
The arc-cost of an arc corresponds to the total objective function cost incurred by group $\group(a)$, and the arc-start-time indicates the time at which group $\group(a)$ depart on a CV.  
These attributes are irrelevant in zero-arcs.  
{A zero-arc indicates when a passenger shares the CV with the next passenger, i.e. is not the last passenger in the group.  This thereby necessitates that any path containing a zero-arc will be followed by a one-arc dictating a group containing the associated passenger.}
\begin{itemize}
	\item Arc-start-time: As mentioned earlier the one-arcs are 
	duplicated for each time instant 
$t \in \cap_{j \in \group(a)} 
[\passengerTimeRequest{j} - \timewindow,\passengerTimeRequest{j} + 
\timewindow]$.
For each such duplicated one-arc 
we have 
$\arcstarttime{a} = t - \timeTerminalToBuilding{d}$,  which 
is a feasible start time from the terminal $\terminal$ 
for the group $\group(a)$. 
    \item Arc-cost: In order to encode the objective function on the arcs, we set 
	\begin{equation}
	\label{eqn:arcCost}
	\arccost{a} := \alpha \cdot \sum\limits_{j \in \group(a)}  \left( \arcstarttime{a}+\timeTerminalToBuilding{\building(j)} - \max\limits_{c \in \trainTrips : \tripTerminalTime{c} \leq \arcstarttime{a}} \trainStartTime{c}{\passengerStop{j}} \right) + (1-\alpha).
	\end{equation}
	The first term is scaled by $\alpha$ and multiplies the total travel time in the objective function, 
which is composed of CV travel time and wating time at $\terminal$. 
\end{itemize}
 
Let $\pathSet^d$ be the set of arc-specified 
$\rootnode^d$-to-$\terminalnode^d$ paths in 
$\DD^d$. For any path $p \in \pathSet^d$, the groups $\groups{p}$ composing the partition 
defined by $p$ is $\groups{p} := \bigcup_{a \in A^d : \arcdomain{a} = 1} \group(a)$.  
	The DDs are constructed in such a way that for every arc-specified 
	$\rootnode^d$-to-$\terminalnode^d$ path, each passenger $j \in \passengerSet(d)$ is in 
	exactly one $g \in \groups{p}$.  	
	The paths also entail the times that each group departs $\terminal$ and the impact on the 
	objective function of selecting an arc.  Time $\arcstarttime{a}$ indicates that the passenger 
	departs to destination $\building$, so that $\group(a)$ occupies the CV assigned to it from time 
	$\arcstarttime{a}$ until $\arcstarttime{a} + \timeRoundTimeToBuilding{\building}$.  We note here 
	that the construction of the DD ensures that the arrival time to destination $d$ for each group 
	is feasible with respect to requested arrival times.  In particular, for each passenger 
	$j \in \group(a)$ we have 
	$\passengerTimeRequest{j} - \timewindow \leq \arcstarttime{a} + \timeTerminalToBuilding{d} \leq \passengerTimeRequest{j} + \timewindow$.

	In the  notation of~(\ref{eqn:IPmodel}), this will correspond to $\sum_{j \in \group(a)} \tottime{j}$ if $z_{j,\arcstarttime{a}} = 1$ for $j \in \group(a)$. The second term in the objective function, $1 - \alpha$, scales the indicator that this represents a CV trip.   The cost of a path $\pathcost{p}$ is the sum of the arc-costs of the one-arcs in $p$.  
	
	Consider the following two properties; $\forall d \in \buildingSet, \forall p \in \pathSet^d, \forall j \in \passengerSet(d)$:
\begin{description}
	\item[\textnormal{\textbf{(DD-1)}}:]  there is exactly one group $g \in \groups{p}$ for which $j \in g$ (denote by $a^p(g)$ the one-arc selecting group $g$);  and
	\item[\textnormal{\textbf{(DD-2)}}:]  for such a group $g \in \groups{p}$ with $j \in g$, $ \passengerTimeRequest{j} - \timewindow \leq \arcstarttime{a^p(g)} + \timeTerminalToBuilding{d}  \leq \passengerTimeRequest{j} + \timewindow$.
\end{description}
We will denote by $g(j) \in \groups{p}$ the unique group to which $j$ belongs in the partition $p$. 
If properties~\textnormal{\textbf{(DD-1)}} and~\textnormal{\textbf{(DD-2)}} are satisfied, then 
any collection $\mathcal{Q}$ of 
$\nbuildings$ \st{$(:= |\buildingSet|)$} paths $\set{p^1, \ldots, p^{\nbuildings}}$, where, for $d \in \buildingSet$, 
$p^d$ is a $\rootnode^d\textnormal{-to-}\terminalnode^d$ in $D^d$, partitions all of $\passengerSet$ 
into $\bigcup_{d \in \buildingSet} \set{\groups{p}}$ groups and the objective function of such a partition 
is $\sum_{d \in \buildingSet} \sum_{p \in \pathSet^d} c(p) = f(\alpha)$. 

Let $\mathcal{G}^d$ be every possible partition of $\passengerSet(d)$ into contiguous subsets for which each subset of passengers can board a common CV.  By Proposition~\ref{prop:optSolutionStructure}, we can consider only these partitions in seeking optimal solutions.  Consider the following property as well:
\begin{description}
	\item[\textnormal{\textbf{(DD-3)}}:] For every partition $\mathsf{g} \in \mathcal{G}^d$, there exists a path $p \in \pathSet^d$ for which $\groups{p} = \mathsf{g}$ and for every time for which the passengers in groups specified by $\groups{p}$ can depart together.
\end{description}
If  property~\textnormal{\textbf{(DD-3)}} is satisfied in $D^d$, then the paths in $\pathSet$ list all possible partitions and departure times from $\terminal$, and therefore defines the feasible region. 

Finally, consider the following property defined over paths 
$\mathcal{Q} = \{p^1, \ldots, p^{\nbuildings}\}$:
\begin{description}
	\item[\textnormal{\textbf{(DD-4)}}:] $\forall t \in \timeSet$, $\left|\bigcup\limits_{d \in \buildingSet} \set{ a \in p^d : \arcstarttime{a} \leq t \leq \arcstarttime{a} + \timeRoundTimeToBuilding{d} } \right| \leq \nCVs $.
\end{description}
If $\mathcal{Q}$ satisfies \textnormal{\textbf{(DD-4)}}, then assigning unique CV trips to each group $g \in \bigcup_{d \in \buildingSet} \set{\groups{d}}$ and leaving $\terminal$ at time $\arcstarttime{g}$ defines a feasible solution to ILMTP-SD that has objective function value $\sum_{d \in \buildingSet} \sum_{p \in \pathSet^d} c(p)$.  Therefore, building a set of DDs satisfying properties~\textnormal{\textbf{(DD-1)}}, \textnormal{\textbf{(DD-2)}} and~\textnormal{\textbf{(DD-3)}}, and finding a collection of paths satisfying condition~\textnormal{\textbf{(DD-4)}} which is of minimum total length outlines another 
approach 
for solving ILMTP-SD.  

Algorithm~\ref{alg:dd} constructs a DD that satisfies properties~\textnormal{\textbf{(DD-1)}} to~\textnormal{\textbf{(DD-3)}} for each destination. 
The algorithm proceeds as follows. 
Algorithm~\ref{alg:dd} starts by computing (Line~\ref{alg:earlylatetimes}) 
for each passenger $j^d_i$, the earliest ($t^{{e}}(j^d_i)$) and latest ($t^{{l}}(j^d_i)$) possible 
departure time 
from $\terminal$ if she were to ride alone using~\eqref{passEarlyLateTimes}.  
Line~\ref{alg:lin:root} creates the root node $\rootnode^d$, which is also referred to as $u_1^0$ for ease of notation. 
Each iteration of the loop in line~\ref{alg:lin:outloop} creates the arcs from layer $i$ to layer 
$i+1$ and the corresponding nodes for every possible state. 
For ease of notation, we also denote a node $u$ in layer $\ell(u) = i$ and state $s(u) = k$ by $u_i^k$.  
Note that for $i = n_d$ the node $u^0_{n_d+1}$ represents the terminal node 
$\terminalnode^d$ of the DD $\DD^d$. 
Line~\ref{alg:lin:0node} creates a single node to add to the $(i+1)$-th layer.  Arcs drawn from the nodes $u^k_i$ in layer 
$i$ to the node $u_{i+1}^0$ are \emph{one-arcs} and represent the grouping of passengers 
$\{j^d_{i-k},\ldots,j^d_k\}$.  The creation of the one-arcs as well as the assignment of CV start times and costs 
are executed in the loop defined by Line~\ref{alg:lin:tloop} of the algorithm.  
The loop in line~\ref{alg:lin:inloop} iterates over the nodes in layer $i$  
and adds \emph{zero-arcs} between the layers $i$ and $i+1$.  
In particular, the loop creates nodes $u^k_{i+1}$ in layer $i+1$ for a positive state $k$, 
which entails grouping passengers $j^d_{i+1-k}, \ldots, j^d_{i+1}$ together through a zero-arc, 
if the conditions in line~\ref{alg:lin:0arc} hold\st{s}: 
(i) the number of passengers does not exceed $\cvCap$; 
(ii) there exists such a passenger $j^d_{i+1}$; 
and (iii) the time windows of passengers $j^d_{i+1-k}$ and $j^d_{i+1}$ overlap.  
Since the passengers are ordered by increasing deadlines, satisfaction of (iii) implies 
that the passengers in $\{j^d_{i+1-k}, \ldots, j^d_{i+1}\}$ have at least one CV starting time 
at $\terminal$ for which they arrive at  their destination within their time windows.    
Note that the cost of arcs is set to zero since this is a zero-arc. Finally, 
the loop in line~\ref{alg:lin:tloop} creates multiple one-arcs for each node, 
which entails that passenger $j^d_i$ is the last in a group,  
by iterating over all possible departure times shared by passengers $j^d_{i-k}$ to $j^d_i$ 
and computing the corresponding cost of such a grouping according to the departure time.

\begin{algorithm}
\footnotesize
\caption{Construction of $\DD^d$ for passengers $j^d_1, \ldots, j^d_{\nPassengersPerBuilding{d}}$ with destination $d$ }
\label{alg:dd}
\begin{algorithmic}[1]
\State Compute for each passenger $j^d_i$ the earliest and latest times of departure from the 
$\terminal$: 
\begin{equation*}\label{passEarlyLateTimes}
\begin{aligned}
& t^{{e}}(j^d_i) \gets \max \set{\passengerTimeRequest{j^d_i} - \timeTerminalToBuilding{d},\min_{c \in \trainTrips} \set{\tripTerminalTime{c}} } \\
& t^{{l}}(j^d_i) \gets \min \set{ \passengerTimeRequest{j^d_i} + \timeTerminalToBuilding{d} - 1 ,\max_{c \in \trainTrips} \set{\tripTerminalTime{c}}}
\end{aligned}
\end{equation*}
\label{alg:earlylatetimes}
\State Add new node $u_1^0$ to $\DDlayer^d_1$ such that $s(u_1^0) = 0$ and $\ell(u_1^0) = 1$
\Comment{Same as root node $\rootnode^d$}
\label{alg:lin:root}
\For{$i \gets 1, \ldots, \nPassengersPerBuilding{d}$}
\label{alg:lin:outloop}
\Comment{Determines transitions after passenger $j_i^d$}
\State Add new node $u_{i+1}^0$ to $\DDlayer^d_{i+1}$ such that $s(u_{i+1}^0) = 0$ and $\ell(u_{i+1}^0) = 1$
\label{alg:lin:0node}
\For{$k \gets 0, \ldots, |\DDlayer^d_i|-1$} 
\Comment{Number of unassigned passengers up to $j_{i}^d$}
\label{alg:lin:inloop}
\If{$k < \cvCap - 1$ \textbf{and} $i < \nPassengersPerBuilding{d}$ \textbf{and} $t^{{e}}(j^d_{i+1}) \leq t^{{l}}(j^d_{i-k})$}
\label{alg:lin:0arc}
\Comment{Checks if $j_{i+1}^d$ can join them}
\State Add new node $u_{i+1}^{k+1}$ to $\DDlayer^d_{i+1}$ such that $s(u_{i+1}^{k+1}) = k+1$ and $\ell(u_{i+1}^{k+1}) = i+1$
\State Add new arc $a$ to $\DDarcs^\building$ such that $\psi(a) = u_i^k, ~\omega(a) = u_{i+1}^{k+1},$ and $\phi(a) = 0$
\State $\eta(a) \gets 0$ 
\Comment{Group cost deferred}
\EndIf
\For{$t \gets t^{{e}}(j^d_{i}), t^{\mathrm{e}}(j^d_{i}) + 1, \ldots, t^{{l}}(j^d_{i-k})$} 
\Comment{Departure times for group $\{ j^d_{i-k}, \ldots, j^d_i \}$}
\label{alg:lin:tloop}
\State Add new arc $a$ to $\DDarcs^\building$ such that $\psi(a) = u_i^j, ~\omega(a) = u_{i+1}^{0}, ~ \phi(a) = 1,$ and $t^0(a) = t$  
\State $\eta(a) \gets  
\alpha \sum\limits_{j \in \group(a)} \left( t + \timeTerminalToBuilding{d} -  \max\limits_{c \in \trainTrips:\tripTerminalTime{c} \leq \arcstarttime{a} }  \set{ \trainStartTime{c}{\passengerStop{j}} } \right) + (1-\alpha)$
\Comment{Group cost incurred}
\EndFor
\EndFor
\EndFor
\end{algorithmic}
\end{algorithm}
 
Theorem~\ref{theorem:DDproperty} below shows that Algorithm~\ref{alg:dd} constructs 
DDs that are polynomial in the size of the input satisfying properties 
\textnormal{\textbf{(DD-1)}} to \textnormal{\textbf{(DD-3)}}.  
Section~\ref{sec:NFR} discusses how to identify the optimal collection of paths satisfying property~\textnormal{\textbf{(DD-4)}}.

 \begin{thm}\label{theorem:DDproperty}
 	For every $d \in \buildingSet$, Algorithm~\ref{alg:dd} constructs DDs $\DD^d$ with $O(\nPassengersPerBuilding{d} \cdot \cvCap)$ nodes and $O(\nPassengersPerBuilding{d} \cdot \cvCap \cdot \timewindow)$ arcs satisfying properties~\textnormal{\textbf{(DD-1)}}, \textnormal{\textbf{(DD-2)}}, and~\textnormal{\textbf{(DD-3)}} that can be constructed in time $O(n \cdot \cvCap \cdot \timewindow)$.
 \end{thm}

The proof of this result is deferred to  Appendix~\ref{app:DDproperty}.


	\begin{example}
		\label{example:1}
		Consider the following ILMTP instance. All passengers request transportation to a single destination $d$. Let $\timeRoundTimeToBuilding{d} = 4, (\timeTerminalToBuilding{d} = 2, \timeStopAtBuilding{d} = 0, \timetBuildingToTerminal{d} = 2)$, $\timewindow = 1$,  $\nCVs = 2$, and $\cvCap = 3$. There are 5 passengers, requesting arrival time to $d$ at 5, 6, 6, 7, 9 for passengers $j^d_1, \ldots, j^d_5$, respectively. Mass transit trips arrive to $\terminal$ at time 2 and 6 ($c=1$ and $c=2$).   All passengers originate from the same $\location$ and 
		$\trainStartTime{c}{\location} = 0, 4$ for $c = 1,2$ respectively.
		
	Figure~\ref{fig:bdd} depicts a DD that satisfies properties~\textnormal{\textbf{(DD-1)}}, \textnormal{\textbf{(DD-2)}}, and \textnormal{\textbf{(DD-3)}}.  The layers are drawn in ascending order from top to bottom.  One-arcs are depicted as solid lines, and zero-arcs as dashed lines, interpreted to be pointing downwards.  Each one-arc \st{$a$ }has two labels depicted in parentheses, the first label is the arc-start-time $\arcstarttime{a}$ and the second label is the total travel time of passengers in $g(a)$.  Each arc-cost (for the one-arcs) is the second argument times $\alpha$ plus ($1-\alpha$).
	
	Consider the red-colored path $p'$, which traverses arcs connecting $\rootnode^b$ to $\terminalnode^b$ through the node-specified path $\rootnode^b-0-0-0-0-\terminalnode^b$ along arcs labeled $(2,4)-(3,5)-(3,5)-(6,4)-(6,4)$.  Each arc emanates from a node with label 0, which indicates that the passengers travel alone in a CV.  The arcs on this path dictates that the passengers leave $\terminal$ at times $t = 2, 3, 3, 6, 6$ and have total travel times of $4,5,5,4,4$ time units, respectively.   To achieve these travel times, passengers $j^d_1, j^d_2, j^d_3$ arrive to $\terminal$ on the mass transit trip $c=1$, and passengers $j^d_4, j^d_5$ arrive on trip $c=2$.  This path satisfies properties~\textnormal{\textbf{(DD-1)}} and~\textnormal{\textbf{(DD-2)}}, as do all paths in the diagram.   However, this path does violate property~\textnormal{\textbf{(DD-4)}}\textemdash consider for example $t=4$. The first three passengers are each assigned CVs that will be en-route at $t=4$ which violates the restriction that $\nCVs = 2$.

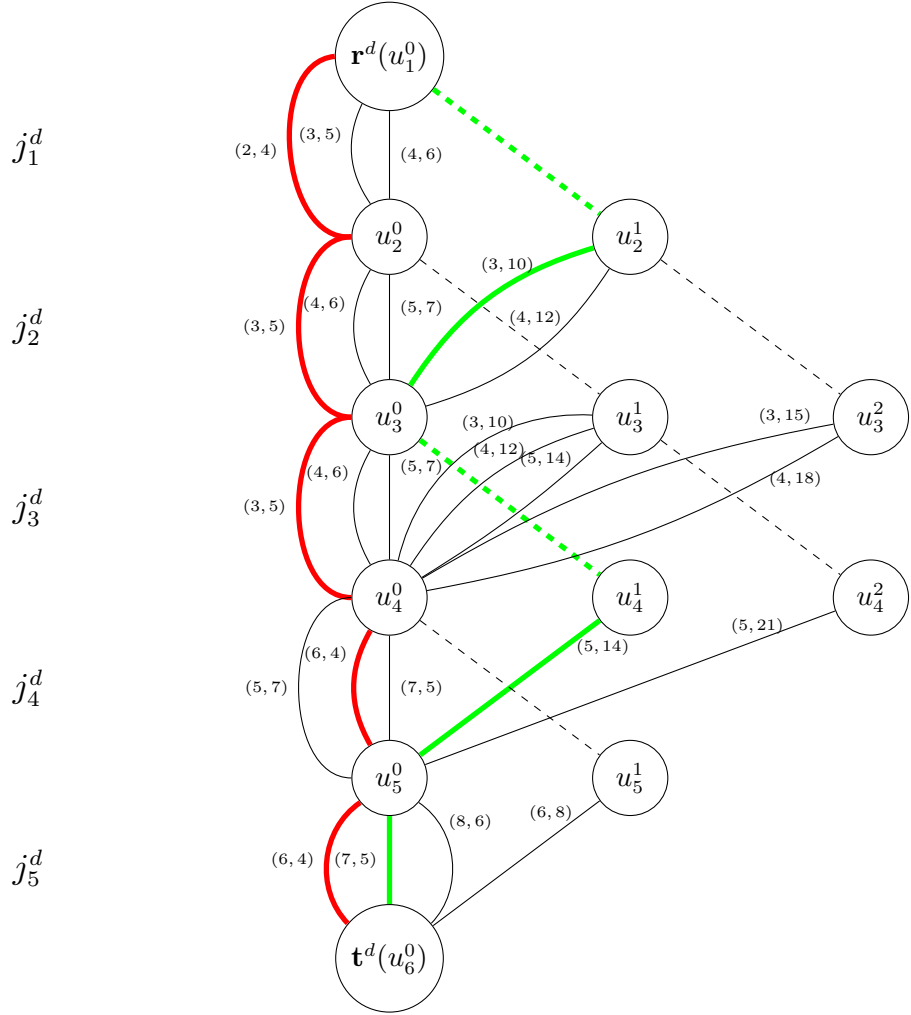
\begin{figure}[t!]
	\centering
	\begin{tikzpicture}[scale=0.35][font=\sffamily]
	
	\node [scale=1.2] (x1) at (-20,-3) {$j^d_1$};
	\node [scale=1.2] (x2) at (-20,-9) {$j^d_2$};
	\node [scale=1.2] (x3) at (-20,-15) {$j^d_3$};
	\node [scale=1.2] (x4) at (-20,-21) {$j^d_4$};
	\node [scale=1.2] (x5) at (-20,-27) {$j^d_5$};
	
	\node [draw,circle,minimum size = 1cm] (er) at (-8,0) {$\rootnode^d (u_1^0)$};
	\node [draw,circle,minimum size = 1cm] (u1) at (-8,-6) {$u_2^0$};
	\node [draw,circle,minimum size = 1cm] (u2) at (0,-6) {$u_2^1$};

	\path[-](er) edge [line width = 2, draw=red, bend right = 90] node [left] {\tiny $(2,4)$} (u1);
	\path[-](er) edge [bend right = 30, pos = 0.3] node [left] {\tiny $(3,5)$} (u1);
	\path[-](er) edge [bend left = 0] node [right] {\tiny $(4,6)$} (u1);
	\path[-](er) edge [line width = 2, draw=green,dashed]  (u2);

	\node [draw,circle,minimum size = 1cm] (u3) at (-8,-12) {$u_3^0$};
	\node [draw,circle,minimum size = 1cm] (u4) at (0,-12) {$u_3^1$};
	\node [draw,circle,minimum size = 1cm] (u5) at (8,-12) {$u_3^2$};
	\path[-](u1) edge [dashed]  (u4);
	\path[-](u2) edge [dashed]  (u5);
	
	\path[-](u1) edge [line width = 2, draw=red, bend right = 90] node [left] {\tiny $(3,5)$} (u3);
	\path[-](u1) edge [bend right = 30, pos = 0.3] node [left] {\tiny $(4,6)$} (u3);
	\path[-](u1) edge [bend left = 0] node [right, pos = 0.3] {\tiny $(5,7)$} (u3);
	
	\path[-](u2) edge [bend right = 20, line width = 2, draw=green] node [left, pos = .2] {\tiny $(3,10)$} (u3);
	\path[-](u2) edge [bend left = 20] node [left, pos = .25] {\tiny $(4,12)$} (u3);

	\node [draw,circle,minimum size = 1cm] (u6) at (-8,-18) {$u_4^0$};
	\node [draw,circle,minimum size = 1cm] (u7) at (0,-18) {$u_4^1$};
	\node [draw,circle,minimum size = 1cm] (u8) at (8,-18) {$u_4^2$};
	
	\path[-](u3) edge [line width = 2, draw=green,dashed]  (u7);
	\path[-](u4) edge [dashed]  (u8);

	\path[-](u3) edge [line width = 2, draw=red, bend right = 90] node [left] {\tiny $(3,5)$} (u6);
	\path[-](u3) edge [bend right = 30, pos = .2] node [left] {\tiny $(4,6)$} (u6);
	\path[-](u3) edge [bend left = 0] node [right, pos = 0.1] {\tiny $(5,7)$} (u6);

	\path[-](u4) edge [bend right = 40, near start] node [left] {\tiny $(3,10)$} (u6);
	\path[-](u4) edge [bend right = 20, near start] node [left] {\tiny $(4,12)$} (u6);
	\path[-](u4) edge [bend left = 5, pos = .1] node [left] {\tiny $(5,14)$} (u6);
	\path[-](u5) edge [bend right = 10,pos=.1] node [above] {\tiny $(3,15)$} (u6);
	\path[-](u5) edge [bend left = 10] node [below, pos = 0.1] {\tiny $(4,18)$} (u6);

	\node [draw,circle,minimum size = 1cm] (u9) at (-8,-24) {$u_5^0$};
	\node [draw,circle,minimum size = 1cm] (u10) at (0,-24) {$u_5^1$};

	\path[-](u6) edge [bend right = 90] node [left] {\tiny $(5,7)$} (u9);
	\path[-](u6) edge [line width = 2, draw=red, bend right = 30] node [left, pos = 0.2] {\tiny $(6,4)$} (u9);
	\path[-](u6) edge node [right, pos = 0.5] {\tiny $(7,5)$} (u9);
	\path[-](u7) edge [line width = 2, draw=green] node [right, pos=.2] {\tiny $(5,14)$} (u9);
	\path[-](u8) edge node [left, pos=.1] {\tiny $(5,21)$} (u9);

	\path[-](u6) edge [dashed]  (u10);
	
	\node [draw,circle,minimum size = 1cm] (t) at (-8,-30) {$\terminalnode^d(u_6^0)$};
	
	\path[-](u9) edge [line width = 2, draw=red, bend right = 50] node [left] {\tiny $(6,4)$} (t);
	\path[-](u9) edge [line width = 2, draw=green]  node [left] {\tiny $(7,5)$} (t);
	\path[-](u9) edge [bend left = 50] node [right, pos = .2] {\tiny $(8,6)$} (t);
	\path[-](u10) edge node [left, pos = .1] {\tiny $(6,8)$} (t);
	
%
%
%
	\end{tikzpicture}
	\caption{Decision diagram for Example~\ref{example:1}} 
	\label{fig:bdd}
\end{figure} 
	
	Consider now the green-colored path $p''$, which traverses arcs connecting $\rootnode^d$ to $\terminalnode^d$ through the node-specified path $\rootnode^d-1-0-1-0-\terminalnode^d$ with one-arcs labeled $(3,10)-(5,14)-(7,5)$ on layers 2, 4, and 5, respectively. This path has three one-arcs that specify groups $\set{j^d_1,j^d_2}, \set{j^d_3, j^d_4}, \set{j^d_5}$ to leave $\terminal$ at times $t = 3, 5, 7$, respectively.  The total travel times on each CV trip are 10 (5+5), 14 (7+7), 5 (5) achieved by passengers $j^d_1, j^d_2, j^d_3, j^d_4$ arriving on mass transit trip $c=1$, and passenger $j^d_5$ arriving on trip $c=2$. For $t = 0, 1, \ldots, 10$, the number of active CVs is 0,0,0,1,1,2,2,2,2,1,0, respectively, upon which all CVs have returned to $\terminal$, thereby never violating the constraints on the number of CVs.   This path therefore satisfies property~\textnormal{\textbf{(DD-4)}} and corresponds to a feasible solution.  The evaluation of the objective function corresponding to this solution depends on $\alpha$, and is evaluated as $\pathcost{p''} = \alpha \cdot(10 + 14 + 5) + (1-\alpha) \cdot 3$.   
	
	There are 492 paths in the depicted DD, corresponding to $|\mathcal{G}^d|$.  This example suggests the advantages of a DD-based approach, in that an exponentially sized set of solutions can be represented, compactly, in a small-sized diagram. 
	
	An additional note is in order. Consider the set of arcs directed between the nodes labeled 0 in the penultimate layers of the DD.  Each arc represents group $\set{j^d_4}$, the singleton passenger traveling alone.  There are three ways this can happen\textemdash the passenger arrives to $\terminal$ on the mass transit trip $c=1$, waits for 3 time units, and then boards a CV, or the passenger arrives to $\terminal$ on trip $c=2$, waits 0 or 1 time units, and then boards a CV.  Should this group be selected as part of the solution, the selection of arc (and, therefore, mass transit trip and waiting time) will depend on the availability of CVs that can be restricted based on CV trips to destination $d$ and to other destinations in the system. 
	\end{example}

%


\subsection{Network-flow reformulation}
\label{sec:NFR}

Given, for each destination $\building \in \buildingSet$, a DD $\DD^d$ satisfying properties~\textnormal{\textbf{(DD-1)}}, ~\textnormal{\textbf{(DD-2)}}, and ~\textnormal{\textbf{(DD-3)}}, we can reformulate the ILMTP-SD optimization problem as a \emph{consistent path problem}.  In each DD, we must select a path so that at any time $t$ no more than $\nCVs$ CVs are assigned to the  groups.  More formally, we want to select, for every $d \in \buildingSet$, a path $p^d \in D^d$ such that, for every $t$, the number of one-arcs with $\arcstarttime{a} \leq t \leq \arcstarttime{a} + \timeRoundTimeToBuilding{\building}$ is less than or equal to $\nCVs$. Let $\arcTimeIndicator{a}{t} \in \set{0,1}$ indicate that a CV would be active at time $t$ (i.e. $\arcstarttime{a} \leq t \leq \arcstarttime{a} + \timeRoundTimeToBuilding{\building}$) if arc $a$ is chosen.  One can formulate this by assigning a variable $y_a$ to each arc $a$ and solving the following optimization problem:

\begin{equation}
\tag{NF}
\label{eqn:spd}
\begin{aligned}
& \min
& & \sum_{\building \in \buildingSet} \sum_{a \in \DDarcs^\building} \arccost{a} y_a \\
& \text{s.t.}
& & \sum_{a : \arcroot{a} = \rootnode^d} y_a = 1,  \quad \sum_{a : \arcterminal{a} = \terminalnode^d} y_a = 1, 
& & \forall d \in \buildingSet	\\
& & & 
\sum_{a : \arcroot{a} = \DDnode} y_a - 
\sum_{a : \arcterminal{a} = \DDnode} y_a  = 0, 
& &
\forall \building \in \buildingSet, 
\forall \DDnode \in \DDlayer^d_2 \cup \ldots \cup  \DDlayer^d_{\nPassengersPerBuilding{d}} \\
& & & 
\sum_{d \in \buildingSet}  \quad \sum_{a \in \DDarcs^d : \arcTimeIndicator{a}{t} = 1} y_a \leq \nCVs, 
& &
\forall t \in \timeSet
\\
& & &
y_a \in \set{0,1} 
& & 
\forall \building \in \buildingSet, \forall a \in \DDarcs^\building
\end{aligned}
\end{equation}

Model~(\ref{eqn:spd}) directly models each DD as a network-flow problem, where we seek to send one unit of flow from $\rootnode^d$ to $\terminalnode^d$.  The sum of the arc weights are minimized, subject to the singular linking constraint, that enforces the restriction on the number of CVs.  Model~(\ref{eqn:spd}) therefore identifies a collection of paths $\mathcal{Q}$ satisfying property~\textnormal{\textbf{(DD-4)}} of minimum total cost, therefore providing a valid formulation for the ILMTP-SD.


\subsection{\fourthrev{A Column Generation Approach}}
\label{sec:bp}

An alternative \secondrev{solution scheme} for the ILMTP-SD given a collection of DDs satisfying properties~\textnormal{\textbf{(DD-1)}}, \textnormal{\textbf{(DD-2)}}, and~\textnormal{\textbf{(DD-3)}} is to use a branch-and-price scheme~\citep{BnP} by associating a binary variable $z_p$ to every $\rootnode^d$-to-$\terminalnode^d$ path in the collection of DDs, where we let $\pathHits{p}{t}$ be the number of one-arcs in $p$ for which $\arcTimeIndicator{a}{t} = 1$ (which we refer to as the \emph{master problem}):

\begin{equation}
\tag{MP}
\label{eqn:mp}
\begin{aligned}
& \min
& & \sum_{\building \in \buildingSet} \sum_{p \in \pathSet^\building} \pathcost{p} z_p \\
& \text{s.t.}
& & \sum_{p \in \pathSet^\building} z_p = 1,
& &
\forall \building \in \buildingSet \\
& & &
\sum_{\building \in \buildingSet} \sum_{p \in \pathSet^d}  \pathHits{p}{t} z_p \leq \nCVs, 
& &
\forall t \in \timeSet
\\
& & & 
z_p \in \set{0,1},
& & 
\forall \building \in \buildingSet, \forall p \in \pathSet^d.
\end{aligned}
\end{equation}  

Since there is an exponential number of variables corresponding to those paths,  we propose to solve this model by \fourthrev{column generation which can be turned into an exact approach via} branch-and-price.  

The procedure begins by defining an initial search-tree node with no branching decisions, and, for all $d \in \buildingSet$, a subset of the paths $\tilde{\pathSet}^d \subseteq \pathSet^d$.  
Let ${\pathSet} = \cup_{d\in\buildingSet}{\pathSet}^d$ and $\tilde{\pathSet} = \cup_{d\in\buildingSet}\tilde{\pathSet}^d$.   
The \emph{restricted master problem} (RMP($\tilde{\pathSet}$)) is~(\ref{eqn:mp}) restricted to only those variables in $\tilde{\pathSet}$.  $\tilde{\pathSet}$ should contain at least one feasible solution, which we address in Section~\ref{sec:initialFeasibleSolution}.

We solve the LP relaxation of~(\ref{eqn:mp}) by \emph{column generation}, where we add paths $p \in \pathSet \backslash \tilde{\pathSet}$ to $\tilde{\pathSet}$ if the associated variable in (\ref{eqn:mp}) has a reduced cost that is negative at the solution corresponding to the optimal LP relaxation of~(RMP($\tilde{\pathSet}$)).  Since we do not  assume an enumeration of $\pathSet$, we identify if such a path 
exists by solving a \emph{pricing problem}~(\ref{eqn:pp}), as described in the following proposition.

\begin{prp}
For all $d \in \buildingSet$, let $\mu_\building$ be the dual variable associating constraint $\sum_{p \in \tilde{\pathSet}^\building} z_p = 1$ and, for all $ t \in \timeSet$, let $\lambda_t$ be the dual variable associated with constraint $\sum_{\building \in \buildingSet} \linebreak \sum_{p \in \tilde{\pathSet}^d}  \pathHits{p}{t} z_p \leq \nCVs$ at an optimal solution to the 
LP relaxation of ~(RMP($\tilde{\pathSet}$)).
For all $d \in \buildingSet, a \in A^b, $ and $t \in \timeSet$, let $\arcTimeIndicator{a}{t}$ indicate if arc $a$ is a one-arc and $t \in \set{\arcstarttime{a}, \ldots, \arcstarttime{a} + \timeRoundTimeToBuilding{d}}$, i.e., that taking the one-arc requires an active CV at time $t$.


Let us define, for all $a \in \cup_{d\in \buildingSet} A^d$, a binary variable $y_a$ and, for all $d \in \buildingSet$, a binary variable $\zeta_\building$.  If the optimal value of 	
\begin{equation}
\tag{PP}
\label{eqn:pp}
\begin{aligned}
& & \min
&  
\sum_{d \in \buildingSet} \sum_{a \in \DDarcs^d} \arccost{a} y_a - \sum_{d \in \buildingSet} \mu_d \zeta_d +  \sum_{d \in \buildingSet} \sum_{a \in \DDarcs^\building} \sum_{t \in \timeSet	} \lambda_t \arcTimeIndicator{a}{t} y_a 
\\
& &
\textnormal{s.t. }
&
\sum_{a : \arcroot{a} = \rootnode^d} y_a = \zeta_d, \quad \sum_{a : \arcterminal{a} = \terminalnode^d} y_a = \zeta_\building, 
& & \forall d \in \buildingSet	\\
& & & 
\sum_{a : \arcroot{a} = \DDnode} y_a - 
\sum_{a : \arcterminal{a} = \DDnode} y_a  = 0, 
& &
\forall \building \in \buildingSet, 
\forall \DDnode \in \DDnodes^d \textnormal{ with } \DDnode \notin \set{\rootnode^\building, \terminalnode^d} \\
& & &
y_a \in \set{0,1} 
& & 
\forall \building \in \buildingSet, \forall a \in \DDarcs^d
\end{aligned}
\end{equation}  
is non-negative, then the optimal LP solution of ~(RMP($\tilde{\pathSet}$)) is an optimal LP solution of~(\ref{eqn:mp}). Otherwise the set of arcs for which $y_a = 1$ defines a path in the 
\secondrev{diagram} $D^d$ for which $\zeta_d = 1$ with negative reduced cost. 
\end{prp}

\proof{Proof.}
Follows immediately from the definition of reduced cost. 
\Halmos
\endproof

Subproblem (\ref{eqn:pp}) decomposes into separate shortest path problems.  For each $\building \in \buildingSet$, let $p^{\building,*}$ be the shortest path and $f^{\building,*}$ be the shortest path length in $\DD^d$, where each arc has length $\arccost{a} - \sum_{t \in \timeSet} \lambda_t \arcTimeIndicator{a}{t}$.   The variable $z_{p^{d,*}}$ associated with the path that achieves the minimum value $f^{\building,*} - \mu_\building$ in the pricing problem will be the variable in the exponential model (\ref{eqn:mp}) that has the lowest reduced cost.  Since this can be done separately for each destination, and since each DD is directed and acyclic, the pricing problem is solved in linear time in the size of the DDs.  Note that the IP solution to~(RMP($\tilde{\pathSet}$)) will always be a feasible solution to the ILMTP-SD instance.  This equips us with a mechanism for generating feasible solution and upper bounds. 

A branch-and-bound search \fourthrev{can be conducted to complete a} branch-and-price \fourthrev{implementation}.  A queue of search-tree nodes $\Gamma$ is defined and initialized as a singleton $\gamma'$.  At any point in the execution of the algorithm, each search node $\gamma \in \Gamma$ is defined by a set of branching decisions $\mathrm{out}(\gamma), \mathrm{in}(\gamma) \in \tilde{\pathSet}$. We also maintain the best-known solution $z^*$ and its objective value $f^*$.  

\fifthrev{While $\Gamma \neq \emptyset$, a search node $\gamma$ is selected to explore. 
The LP relaxation of~(RMP($\tilde{\pathSet}$)), with additional  constraints requiring $z_p = 0,1$ for those paths $p \in \mathrm{out}(\gamma),\mathrm{in}(\gamma)$, respectively, is solved via column generation. This can be implemented as follows.  When $z_p=1$, the associated DD can be simplified to a single path.  When $z_p=0$, the associated DD can be modified to delete only the path $p$ and retain all other paths, noting that removing one path from a DD takes linear time and increases the width of the DD by at most one \citep{CirHoo2014}. One can store the collection of DDs resulting from the branching decisions at each search tree node, or remove the paths with $z_p=0$ each time, starting from the original DDs. 
Suppose $\{p^{\building,1},\ldots,p^{\building,l}\} \subset \mathrm{out}(\gamma)$ be an arbitrary ordering of all paths in $\mathrm{out}(\gamma)$ that are in $\DD^\building$ and are to be removed, i.e. $z_p = 0$ for all $p \in \{p^{\building,1},\ldots,p^{\building,l}\}$. Note that $l$ can vary by $\building$. The procedure for removing the said paths from $\DD^\building$ is iterative.  Starting with $\DD^{\building,0} = \DD^{\building}$, at each iteration $k = 1,\ldots,l$ the path $p^{\building,k}$ is removed from $\DD^{\building,k-1}$ as described by \citep{CirHoo2014} and the resulting diagram is denoted as $\DD^{\building,k}$.  The diagram $\DD^{\building,l}$ is the decision diagram for $\building$ that is used in the relaxation of node $\gamma$. 
If the optimal value of the LP relaxation of~(RMP($\tilde{\pathSet}$)) is greater than or equal to  $f^*$, the node is pruned, and search continues by selecting another node in $\Gamma$.  Otherwise, the IP~(RMP($\tilde{\pathSet}$)) is solved and if the optimal value $f'$ is less than $f^*$, this solution replaces $z^*$ and $f^*$ is updated with $f'$. Note that enforcing a path is the same as removing a DD from the problem and fixing this decision, and since there is a one-to-one mapping between passenger-groups traveling to a building and paths in a DD, this will not cycle and is complete. We note here that our experiments did not require branching since the root node gap is shown experimentally to be very small. }


We also describe another approach to 
identifying a feasible solution in Section~\ref{sec:constructFeasible} since the solution of the   
IP~(RMP($\tilde{\pathSet}$)) can be computationally prohibitive. \secondrev{In general, a branching scheme must be devised in order to close any gap at the root node.  As discussed in the experimental results, the optimality gap at the root node is very small in this application, and so we need not explore branching rules. } 

\subsection{Finding an initial feasible solution}
\label{sec:initialFeasibleSolution}

We can generate an initial feasible solution to~(\ref{eqn:mp}) by appending an extra path to each DD that represents not assigning any passenger.   Specifically, for each $d \in \buildingSet$, create nodes $\DDnode^d_2, \ldots, \DDnode^d_i$ with state $\DDstate{\DDnode^d_i} = \emptyset, \forall i =2, \ldots, n$.   Add one-arcs from $\rootnode^d$ to $\DDnode^d_2$, from $\DDnode^d_n$ to $\terminalnode^d$, and, for $i=2, \ldots, n-1$, from $\DDnode^d_i$ to $\DDnode^d_{n-1}$.  For each new arc $a$ set $\arccost{a} = -\infty$ and $\arcstarttime{a} = 0$.  For the new path $p^d_0$ in each DD set $z_{p^d_0} = 1$ and all other variables equal to 0.  Since for all $t\in \timeSet, \pathHits{p^d_0}{t} = 0$, 
this will be an initial feasible solution to~(\ref{eqn:mp}).  

\subsection{Identifying a feasible solution}
\label{sec:constructFeasible}

Suppose that 
$z^*$ is a solution to the LP relaxation of~(RMP($\tilde{\pathSet}$)) and that 
$z^*$ is not integral.  Denote by $p^d \in \tilde{\pathSet}^d$ for each destination $d \in \buildingSet$ 
as the path satisfying $p^d = \arg\max_{p \in \tilde{\pathSet}^d} z^*_p$.  In other words, 
$p^d$ is path for destination $d$ with the largest fractional value in the solution of the LP relaxation.  
As defined in Section~\ref{sec:singleDestBDD}, let $\groups{p^d}$ represent the partition 
of $\passengerSet(d)$ defined by $p^d$.  The path $p^d$ encodes the particular start times 
on the CV for the groups in $\groups{p^d}$.  In the following, we describe \secondrev{an} IP that fixes the 
groups in $\groups{p^d}$ but attempts to assign starting times for the groups such that the resulting 
CV trips are feasible for the ILMTP-SD.  Note that by assigning different start times we 
are implicitly enumerating other paths in the DD $\DD^d$ with the same groupings. 
The IP formulation is a simplified version of the one presented 
in~\cite{iflmpicaps18}, where route assignments were also considered in the IP.

Prior to describing the model we introduce relevant notation.  Define earliest and latest 
CV start times for the group $\group \in \groups{p^d}$, for all $d \in \buildingSet$, as 
$t^{{e}}(\group) = \max_{j \in \group} t^{{e}}(j)$ and 
$t^{{l}}(\group) = \min_{j \in \group} t^{{l}}(j)$ where $t^{{e}}(j), t^{{l}}(j)$ 
are as defined in~\eqref{passEarlyLateTimes}.  Define the objective value  
associated with the particular start time $t \in [t^{{e}}(\group),t^{{l}}(\group)]$ as 
\[
\eta(\group,t) = \alpha \sum\limits_{j \in \group} \left( t+\timeTerminalToBuilding{d} - 
\max\limits_{c \in \trainTrips :  \set{\tripTerminalTime{c}} \leq t}  \set{\trainStartTime{c}{\passengerStop{j}}} \right) + (1-\alpha).
\]
Further, let $\chi(\group,t,t') \in \set{0,1}$ indicate the times $t'$ (i.e.  $t \leq t' \leq t+
\timeRoundTimeToBuilding{\building}$) at which a CV serving group $\group$ leaving terminal at 
time $t$ would be active.  
The decision variable in the IP formulation is $x_{\group,t} \in \{0,1\}$ for 
$t \in  [t^{{e}}(\group),t^{{l}}(\group)]$ indicating the choice of start time of $t$ on the CV for 
group $\group$. 
The IP formulation is:
{\begin{subequations}\label{optallcv.form}
\begin{align}
	\min
&\,	 \sum\limits_{\group \in \bigcup\limits_{d \in \buildingSet} \groups{d}} \eta(\group,t) x_{\group,t} \\
	\text{s.t.}							&\,	\sum\limits_{t = t^{\mathrm{e}}(\group)}^{t^{\mathrm{l}}(\group)}  
										x_{\group,t} = 1\, && \forall\, \group \in \bigcup\limits_{d \in \buildingSet} \groups{d}
										\label{oneronet} \\
									&\,	\sum\limits_{\group \in \cup_{d \in \buildingSet} \groups{d}} 
									\sum\limits_{t =  t^{{e}}(\group)}^{t^{{l}}(\group)} \chi_{\group,t,t'} x_{\group,t}
												\leq \nCVs \, &&\forall\, t' \in \timeSet.		\label{restrictNumCVs}
\end{align}   
\end{subequations}
Constraint~\eqref{oneronet} imposes that each group is assigned to exactly one CV route and start time.  
Constraint~\eqref{restrictNumCVs} ensures that the number simultaneous trips on the 
CVs does not exceed the number of CVs.  
A feasible solution to~\eqref{optallcv.form} is easily seen to be a feasible solution ILMTP-SD.  
The above IP typically solves at the root node in fractions of a second.


\section{Numerical Evaluation}
\label{sec:numericalEval}

This section details an experimental evaluation conducted to evaluate the efficiency of the algorithms presented, sensitivity analysis on the quality of the solutions obtained, and how the solutions to the ILMTD-SD can be used in more general settings. Section~\ref{sec:instanceGen} describes the instances\thirdrev{, all of which are available at \url{ftp://ftp.merl.com/pub/raghunathan/LastMile-TestSet/}.}  Due to company policy, we are unable to make the data and code available at the INFORMS Journal of Computing repository, but we do provide a sample instance with the structure as those used in the experiments (\url{https://github.com/INFORMSJoC}). Section~\ref{sec:algocomp} provides 
an 
experimental evaluation of the formulations and algorithms. 
We analyze the trade-offs on different objectives and the effect of the time window length on the quality of the solutions obtained in Sections~\ref{sec:obj} and \ref{sec:tw}, respectively. 
Finally, 
in Section~\ref{sec:progen}
we study the quality of the solution to the ILMTP-SD for more general variants of the problem.
All experiments were run on a machine
with an Intel(R) Core(TM) i7-4770 CPU @ 3.40GHz and
16 GB RAM. All algorithms were implemented in \texttt{Python
2.7.6} and the ILPs are solved using \texttt{Gurobi 8.0.1}.

\subsection{Instance generation}
\label{sec:instanceGen}

We generate instances with number of destinations $\nbuildings \in \set{10,25,50}$.  In order to test how well the algorithms scale, we specify the number of passengers 
per destination as $\frac{\nPassengers}{\nbuildings} \in \set{100,150,200}$, and use the corresponding value for $\nPassengers$ (i.e., if $\nbuildings =10$ and $\frac{\nPassengers}{\nbuildings} = 100$, we use $\nPassengers = 1000$) so that our instances have up to 10,000 passengers. We set the number of CVs $\nCVs = \mathsf{round}(0.06*n)$ where $\mathsf{round}(\cdot)$ rounds to the nearest integer and CV capacity $\cvCap = 5$.  We generate 5 instances per configuration.  The number of 
stations where passengers board the mass transportation system is $4$ and so $\stationSet = \{\terminal,1,2,3,4\}$.  The  station of origin for each passenger is 
generated  independently and uniformly at random from $\{1,2,3,4\}$. 
We discretize time in 2 units per minute and we consider serving passengers with requests within a full hour period, 
hence totalling 120 time units of requests, 
which span between 90 and 210. 
Hence, the requested arrival time
$\passengerTimeRequest{\cdot}$ is generated independently and uniformly at random from the set $\set{90,91, \ldots, 210}$.  
The trips $\trainTrips$ consist 
of trains that depart the farthest station $4$ at times $\trainStartTime{c}{4} \in \{0,30,60,\ldots,210\}$.  
The travel time between stations is 10 and so $\trainStartTime{c}{s} = \trainStartTime{c}{4} + (4-s) \cdot 10$ 
for $s = 1,2,3$ and $\tripTerminalTime{c} = \trainStartTime{c}{4} +40$ for all $c \in \trainTrips$.  
The travel time from $\terminal$ to destinations $d$ is chosen independently and uniformly at 
random between $\set{10, 11, \ldots, 20}$, so if travel time to a destination is $t^d$ then 
$\timeTerminalToBuilding{d} = t^d +1$ where $1$ is the time to board, 
$\timeStopAtBuilding{d} = 1$, which is the time to deboard, and 
$\timetBuildingToTerminal{d} = t^d$.

\secondrev{We will at certain points of the discussion evaluate the impact of choosing different values of $\alpha \in [0,1]$. Since the two objectives are on a different scale, we multiple the coefficient of the number of trips by 100 in order to put the objective on a relatively equivalent scale.}


\subsection{Algorithm comparison}\label{sec:algocomp}

We test the efficiency of our three algorithms.  \textbf{IP} refers to solving~\eqref{eqn:IPmodel}.  \textbf{NF} refers to solving~\eqref{eqn:spd}.  \textbf{CG}  refers to solving~(\ref{eqn:mp}) through the column generation algorithm in Section~\ref{sec:bp}. For comparison with heuristic approaches, \textbf{BAS} refers to the break-and-shift heuristic introduced by~\cite{iflmpicaps18}. \textbf{BAS} first sorts the passengers for each destination by their arrival times, creates groups sequentially from such ordering by adding passengers to each new group while possible, then iteratively breaks groups in which one passenger delays another and move passengers to consecutive groups when no travel time delay is incurred. The four algorithms are applied to each generated instance and for each setting of $\timewindow$ and $\alpha$, resulting in 450 runs each.  Since the transportation systems will in practice be repeatedly optimized, the ILMTP-SD requires efficient solution 
methodologies.  Accordingly, we set a time limit of 10 minutes. 

\textbf{IP} does not solve any of the instances tested to optimality within 10 minutes of computational time, and only finds a feasible solution to \secondrev{2} instance\secondrev{s within the time limit without closing the optimality gap.  Even when supplied with a warm-start solution provided by the heuristic from \cite{iflmpicaps18} and adding ordering constraints based on Theorem~\ref{prop:optSolutionStructure}, \textbf{IP} fails to solve any instance tested to optimality within 10 minutes}.  All other algorithms solve all instances in 10 minutes.  
This provides clear indication of the superiority of the DD-based algorithm, and so for the remainder of this section we provide a comparison only of \textbf{NF} and \textbf{CG}.    We note that the root-node optimality gap ($\frac{UB-LB}{LB} \times 100$ where $LB$ is the lower bound and $UB$ is the upper bound) for the DD-based model is below 0.5\%, 
demonstrating the quality of the LP relaxation.  We therefore only report results for solving the root node as branching is not necessary, but can be used if closing the gap is critical.  

Figure~\ref{fig:NFVersusBP-PP} depicts a cumulative distribution plot of performance over all runs for \textbf{NF}, \textbf{CG}, and \textbf{BAS}.    Each line corresponds to an algorithm, and each point on a line is composed of coordinates $(t,s)$ which is the number of instances solved $s$ by time $t$ by the algorithm.  The figure demonstrates that \textbf{CG} solves more problems than \textbf{NF} for any time limit. 
In addition, \textbf{CG} finds some optimal solutions faster than the \textbf{BAS} runtime.
%

A more detailed pairwise comparison of \textbf{CG} with \textbf{NF} and \textbf{BAS} appears in the plots of Figure~\ref{fig:NFVersusBPandBAS}.   The coordinates are the runtime of \textbf{NF} and the runtime of \textbf{CG} and \textbf{BAS}, respectively, in the left and right plot, in seconds in log-scale.  The size of each point correspond to $\nPassengers$ (increasing in size as $\nPassengers$ increases) \secondrev{and} the color of each point corresponds to $\alpha$. 
This plot more readily reveals the advantage of \textbf{CG}.  In only 17 of the smallest instances is the the runtime of \textbf{NF} smaller than that of \textbf{CG}.  Additionally, the relative superiority of \textbf{CG} grows as $\alpha$ decreases and the objective is more scaled to the number of CV trips.  A similar comparison of runtime exists for \textbf{CG} versus \textbf{BAS}, except that for smaller instances \textbf{BAS} is more efficient, while for larger instances \textbf{CG}, an exact approach, can 
be more efficient than the heuristic \textbf{BAS}. 

\begin{figure}[t!] 
\centering
		\includegraphics[scale=0.05]{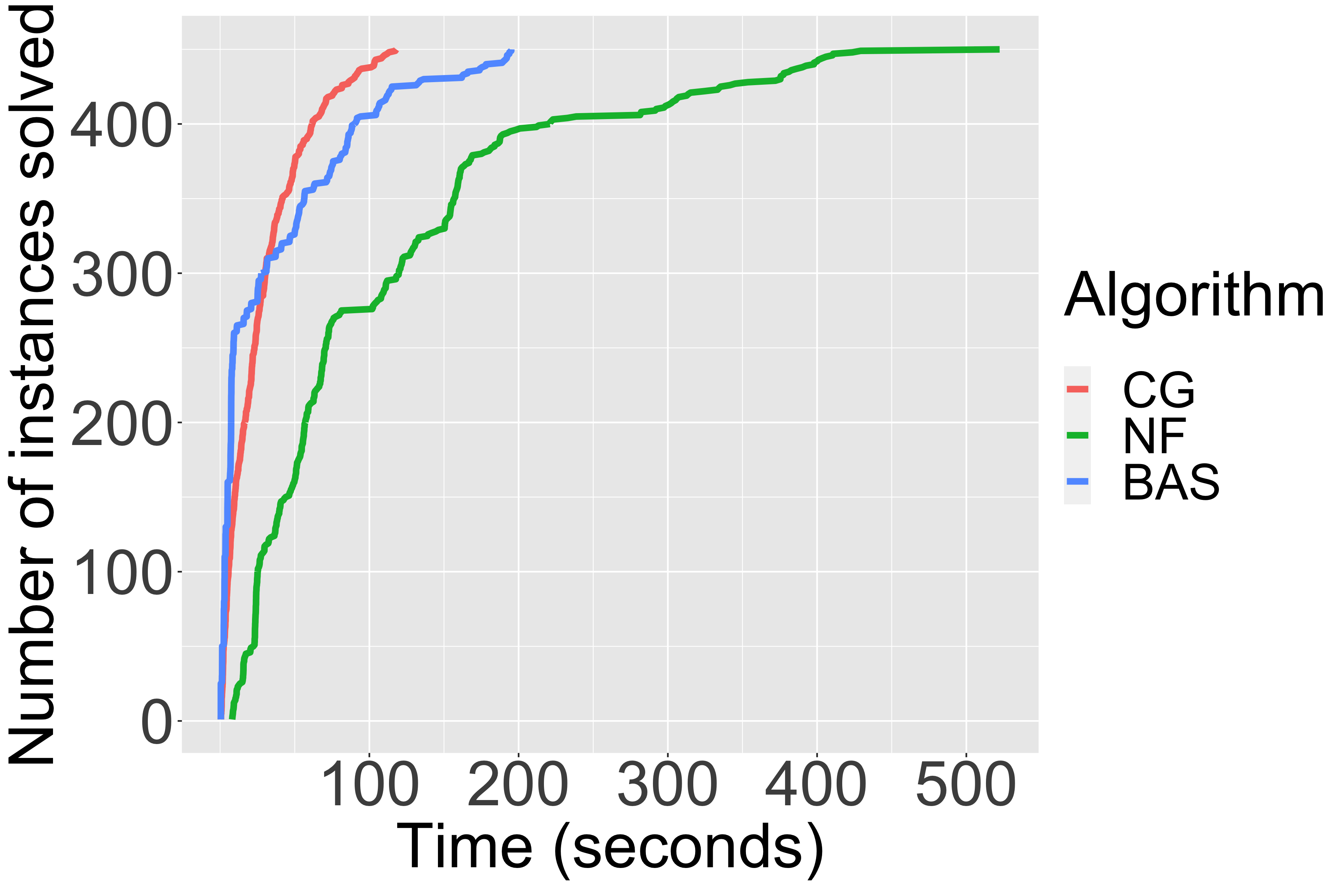}
	\caption{Cumulative distribution plot of performance, comparing \texttt{CG} and \texttt{NF}.} 
	\label{fig:NFVersusBP-PP}
\end{figure}

 
\begin{figure}[t!] 
	\begin{minipage}[b]{0.5\linewidth}
		\centering
		\includegraphics[scale=0.08]{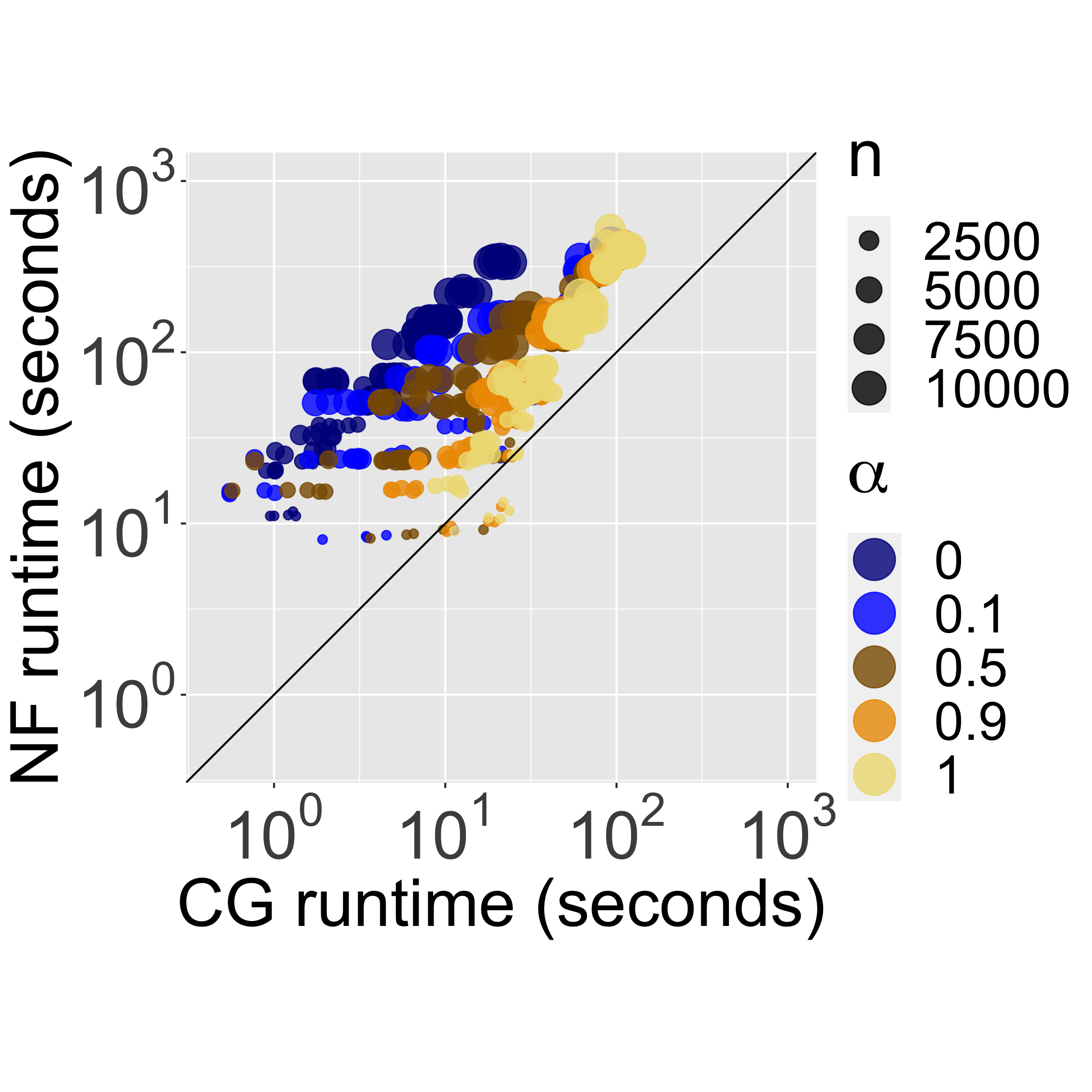}
	\end{minipage}
	\begin{minipage}[b]{0.5\linewidth}
		\centering
		\includegraphics[scale=0.08]{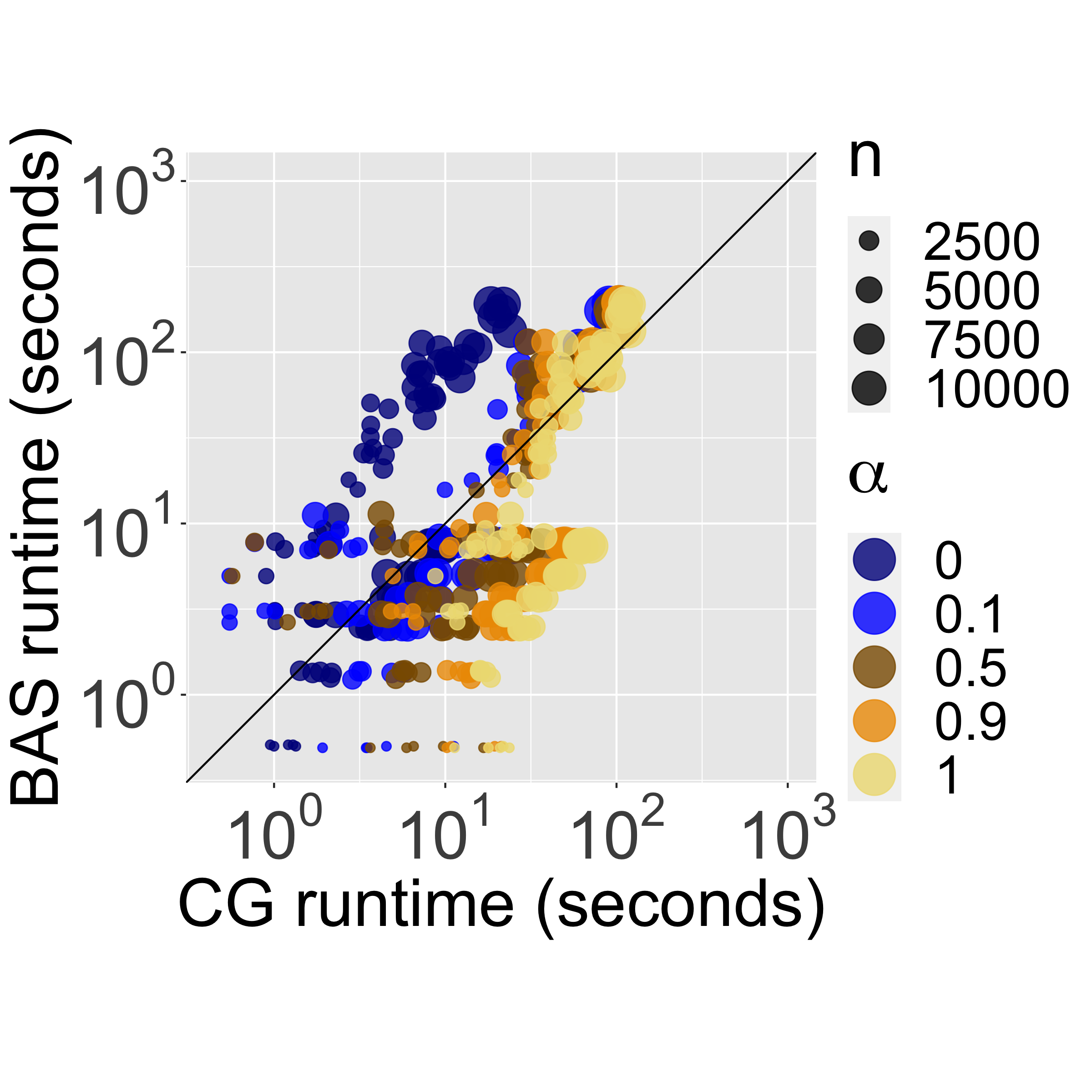} 
	\end{minipage}	
	\caption{Scatter plots, comparing \textbf{CG} with \textbf{NF} (left) and \textbf{CG} with \textbf{BAS} (right).} 
	\label{fig:NFVersusBPandBAS}
\end{figure}

It is also noteworthy  to compare	
the relative quality of the solutions obtained by \textbf{BAS} to the exact approach \textbf{CG}. 
The plot in Figure~\ref{fig:NFVersusBP-Scatter} depicts a scatter plot with the same coordinates as those in the right plot in Figure~\ref{fig:NFVersusBP-Scatter} comparing the runtimes of \textbf{CG} and \textbf{BAS}, but this time colored by the percent improvement in objective function value of the solution obtained by \textbf{CG} over those obtained by \texttt{BAS}, which is always non-negative because \textbf{CG} is an exact approach \secondrev{ and always finds optimal solutions with the 10 minute timeout}. The interesting observation is that for small $\alpha$, \textbf{CG} and \textbf{BAS} obtain nearly the same solution quality, but \textbf{CG} is more efficient.  For large $\alpha$, \textbf{CG} takes slightly longer than \textbf{BAS} but obtains higher quality solutions.  We therefore observe that for any setting of $\alpha$ there are gains that can be realized by utilizing the exact approach. 

%
%


\begin{figure}[t!] 
		\centering	\includegraphics[scale=0.3]{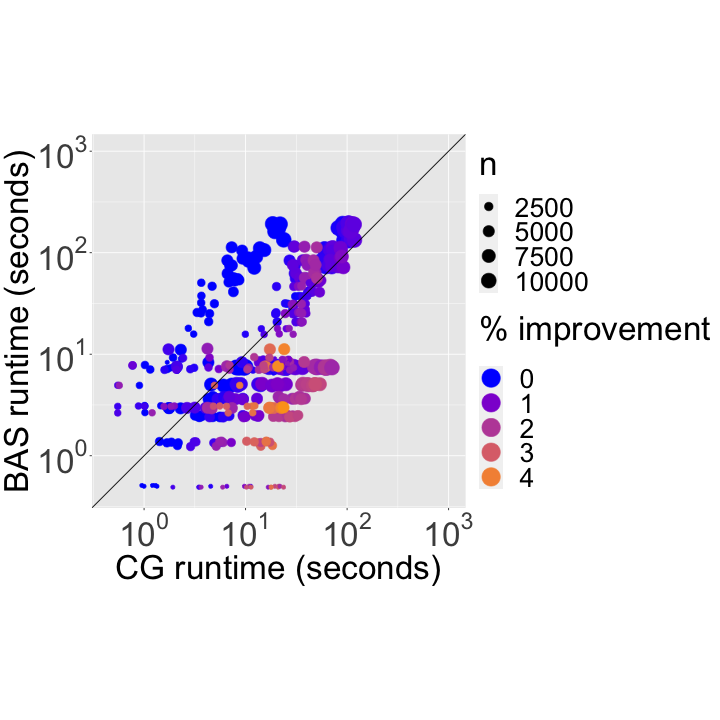}
	\caption{Scatter plot of comparing \textbf{CG} and \textbf{BAS}, but colored by percent improvement in solution quality obtained by \texttt{CG} over that obtained by \texttt{BAS}. } 
	\label{fig:NFVersusBP-Scatter}
\end{figure}

\begin{figure}[t!] 
		\centering
	\includegraphics[scale=0.3]{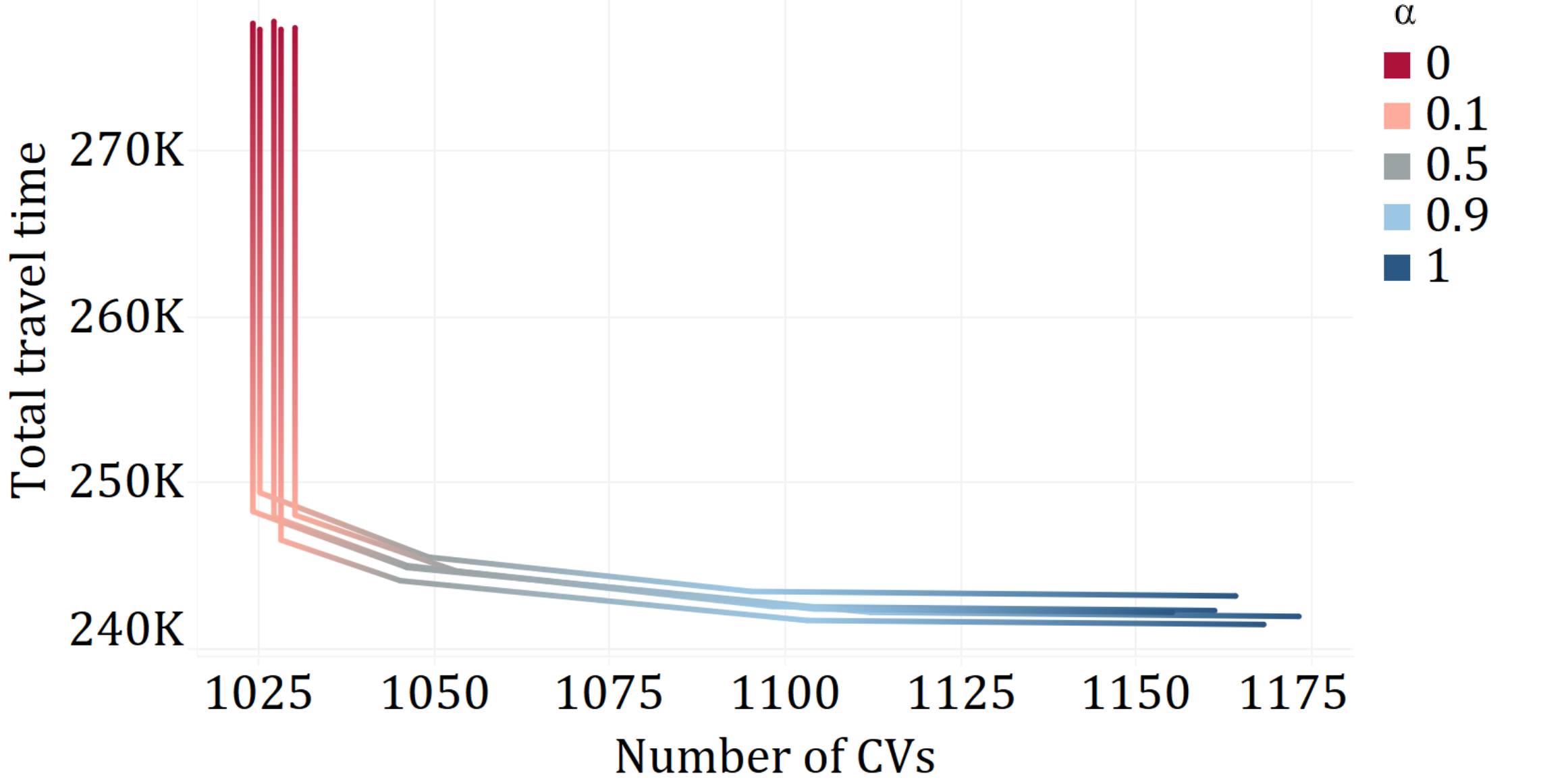} 
	\caption{ Line plot depicting total passenger travel time and number of CVs for different settings of $\alpha$.} 
	\label{fig:ObjectiveWeights}
\end{figure}



\subsection{Objective trade-off analysis}\label{sec:obj}

It is critical to understand how the two objectives, total passenger wait time and number of CV trips, affect the solutions obtained, so that proper operation can be determined.  As \textbf{CG} is shown to be the best algorithm among those tested, we use solutions obtained by \textbf{CG} for the remainder of this section. The \secondrev{plot in} Figure~\ref{fig:ObjectiveWeights} depicts the data through a line plot. Each point corresponds to the total passenger travel time and the number of CVs trips in the solutions obtained by \textbf{CG}, for those instances with $\timewindow = 5, \nPassengers = 5000$ and $\nbuildings = 50$, where we scaled the number of 
trips by $100$ to ensure that the values of the two objectives are of the same order of magnitude. The color corresponds to $\alpha$.


This plot highlights the power of considering a balanced objective.  First, when $\alpha \in \set{0,1}$, we see that considering only one objective can result in great solution for that metric alone, but can be very bad for the other, ignored objective.  Changing $\alpha$ only slightly away from the boundary (i.e., to $\alpha \in \{0.1,0.9\}$) sacrifices only a little on the main objective. 

 The plot reveals that setting $\alpha = 0.1$ leads, in general, to the most balanced solutions and so operators of systems should consider weighting the objectives in this region.  In particular, over all instances tested, there is no difference in the number of CV trips in the optimal solutions obtained when changing $\alpha = 0.0$ to $\alpha = 0.1$.  The same change of $\alpha$ results in a reduction of total travel time from 236,693 second to 199,341, on average over all instances, a decrease of 15\%.  We also see a significant drop in number of CV trips as we change alpha from 1.0 to 0.9.  On average, the number of CV trips decreases from 928 to 862 (a 7.1\% decrease), while only resulting in an increase of average total travel time from 197,273 to 199,341 (a 1.0\% increase).   This will therefore lead to a significant decrease in operations costs and environmental impact with only slightly longer total passenger travel time, and so should be employed in operational decision making.

\subsection{Effect of time-window variation on solutions}\label{sec:tw}

This subsection provides an analysis of how varying the time window affects the solutions obtained.  In particular, we consider the solution time and quality of solutions obtained over all instances tested, for $\timewindow = 5$ and $\timewindow = 10$.
%
%
%
In general, the quality of the solutions obtained are only marginally different when we widen the time windows.  Averaged over all instances, the number of CV trips is 888 and 877 and the total travel time is 213,018 and 198,535, for $\timewindow = 5$ and $\timewindow = 10$, respectively.  This represents a decrease in  number of CV trips of 1.21\% and a decrease in total travel time of 6.80\%.  The increase in solution time is much more substantial, increasing from 14.44 to 40.05 seconds, on average a 177.41\% increase.  This indicates that allowing more flexibility in arrival time constraints makes the problem significantly harder, but results in slightly better operational decisions and so an operator might try to solve the problem with relatively large time windows, but then decrease this flexibility should solutions need to be obtained expediently. 

\subsection{\secondrev{Effect of increasing CV capacity}} \label{sec:increaseCap}

\secondrev{We also evaluated the solution quality and the solution time of \textbf{CG} on instances with increasing vehicle capacity. We selected instances with medium sized attributes, namely $\nbuildings = 25, \timewindow = 10, \nPassengers = 3,750$ (so that $\frac{\nPassengers}{\nbuildings} = 150$), and hence $\nCVs = 225$, using the 5 generated instances per configuration.  We tested with $\cvCap \in \set{5, 6, 7, 8, 9, 10}$.}

\secondrev{The results are summarized in Table~\ref{tab:increaseCap}.  For a fixed $\alpha$, the solution times are generally longer with lower vehicle capacity.  The solution changes significantly as $\cvCap$ grows for small $\alpha$, whereas the differences in the optimal solution for different $\cvCap$ are marginal when $\alpha$ is larger.  This can be attributed to the fact that, as $\alpha$ grows, the emphasis in the objective shifts to total wait time, and CVs will depart the terminal as quickly as possible --- even if that means not filing capacity. Note that the problems are harder to solve at lower capacity.}

\begin{table}[]
\centering
\footnotesize
\begin{tabular}{>{\color{black}}>{\color{black}}l>{\color{black}}l|>{\color{black}}c>{\color{black}}c>{\color{black}}c}
$\alpha$ & $\cvCap$ & Average travel time & Average number of CV trips & Average solution time \\
\hline
0     & 5   & 209,030.2           & 759.2           & 3.724             \\
0     & 6   & 208,846.0             & 635.6           & 1.706             \\
0     & 7   & 208,928.6           & 547.6           & 1.986             \\
0     & 8   & 209,162.0             & 481.4           & 2.178             \\
0     & 9   & 209,575.4           & 431.2           & 2.624             \\
0     & 10  & 209,864.6           & 392.4           & 2.796             \\
\hline
0.1   & 5   & 186,014.8           & 759.2           & 6.426             \\
0.1   & 6   & 186,897.8           & 635.6           & 1.772             \\
0.1   & 7   & 188,177.8           & 547.6           & 2.034             \\
0.1   & 8   & 189,933.8           & 481.4           & 2.392             \\
0.1   & 9   & 191,461.0             & 431.2           & 2.514             \\
0.1   & 10  & 192,911.2           & 392.4           & 2.626             \\
\hline
0.5   & 5   & 184,019.6           & 769.6           & 11.664            \\
0.5   & 6   & 184,155.4           & 651.6           & 1.814             \\
0.5   & 7   & 184,815.8           & 565.6           & 1.950              \\
0.5   & 8   & 185,433.4           & 503.4           & 2.128             \\
0.5   & 9   & 185,808.2           & 459.4           & 2.312             \\
0.5   & 10  & 186,391.4           & 421.0             & 2.584             \\
\hline
0.9   & 5   & 183,060.8           & 792.4           & 22.826            \\
0.9   & 6   & 181,419.2           & 742.8           & 3.668             \\
0.9   & 7   & 181,665.6           & 671.4           & 1.886             \\
0.9   & 8   & 181,720.4           & 638.4           & 2.210              \\
0.9   & 9   & 181,817.4           & 610.8           & 2.314             \\
0.9   & 10  & 181,861.4           & 591.6           & 2.590              \\
\hline
1     & 5   & 182,952.0             & 822.2           & 29.004            \\
1     & 6   & 181,029.0             & 819             & 18.448            \\
1     & 7   & 180,719.2           & 821.4           & 18.182            \\
1     & 8   & 180,573.4           & 822.4           & 18.942            \\
1     & 9   & 180,491.6           & 823.4           & 20.844            \\
1     & 10  & 180,430.0             & 822.8           & 21.520            
\end{tabular}
\caption{\secondrev{Comparison of solution quality and runtime of \textbf{CG} for increasing CV capacity.}}
\label{tab:increaseCap}
\end{table}

\subsection{Problem generalizations}\label{sec:progen}

There are various assumptions we place on the transportation system in order to ensure the optimality of the solutions obtained by \textbf{CG}.  In particular, we require the structural results from \S~\ref{sec:solnStructure}.  Although the proof of optimality remains valid when we extend to time-dependent travel times, there are other dimensions of the problem that we can expand which results in  \textbf{CG} returning only heuristic solutions.  As discussed in~\cite{iflmpicaps18}, the solutions obtained can still be close to optimal even in more general settings.  This subsection provides an analysis of two main problem generalization, and uncovers that the solutions we identify through \textbf{CG} are often superior to what can be found by other techniques, even in the more general settings where the solutions may not be optimal. 

\subsubsection{Express trains} The structural result from \S~\ref{sec:solnStructure} which proves that there exists an optimal solution where passengers going to a common destination are partitioned in the order of their arrival times fails to remain true when there are more complex train systems that are integrated with last-mile transportation.  For example, if trains do not arrive sequentially at regular intervals, stopping at each and every station, this result no longer holds. This is a common characteristic of real-world train systems, where express trains increase the frequency that high-traffic train stops have train service. 

In order to test how well \textbf{CG} performs as a heuristic in this setting, we generate instances with express trains as follows.  We generated instances with $\nbuildings \in \{10,25\}$  and specify the 
number of passengers per destination $\frac{\nPassengers}{\nbuildings} \in \{50,75\}$.  We set $\nCVs = \mathsf{round}(0.1*\nPassengers)$ 
which is a larger number of CVs as a fraction of passengers than the used in the previous tests.  
This was done primarily to provide the \text{MIP} formulation the benefit of solving more of 
the instances.  We generate 5 instances per configuration.  The  set of passengers 
are generated independently at random as described in the beginning of the section.  The trips for the trains consist of express trips that stop only at stations 4 and 2 before reaching the $\terminal$.  These trips occur every 30 time units and can be specified 
using the notation of the paper as: 
$\trainStartTime{c}{4} \in \{0,30,\ldots,180\}$, $\trainStartTime{c}{2} = \trainStartTime{c}{4}+10$,  
$\tripTerminalTime{c} = \trainStartTime{c}{4}+20$ and $\trainStartTime{c}{s} = -\infty$ for $s=1,3$.  
Another set of express trips stop only at stations 3 and 1 before reaching the $\terminal$.  
These trips also occur every 30 time units and can be specified as:
$\trainStartTime{c}{4} \in \{20,50,\ldots,200\}$, $\trainStartTime{c}{1} = \trainStartTime{c}{3}+10$,  
$\tripTerminalTime{c} = \trainStartTime{c}{4}+15$ and $\trainStartTime{c}{s} = -\infty$ for $s=2,4$.   We refer to \textbf{IP} as the using an exact model for this problem setting \citep{iflmpicaps18}.

\begin{table}[]
\centering
\footnotesize
\begin{tabular}{lll|lll|lll}
$\nPassengers$ & $\nbuildings$ & $\alpha$ & $n_{\textbf{CG}}$ & $\hat{t}_{\textbf{CG}}$ & $\hat{z}_{\textbf{CG}}$ & $n_{\textbf{IP}}$ & $\hat{t}_{\textbf{IP}}$ & $\hat{z}_{\textbf{IP}}$ \\ \hline
500            & 10            & 0        & 5                 & 0.1                     & 11760.0                 & 0                 & 601.5                   & 11760.0                 \\
500            & 10            & 1        & 5                 & 1.6                     & 18691.4                 & 5                 & 33.7                    & 18535.2                 \\
750            & 10            & 0        & 5                 & 0.2                     & 15780.0                 & 0                 & 603.8                   & 15880.0                 \\
750            & 10            & 1        & 5                 & 2.4                     & 27095.8                 & 5                 & 14.5                    & 26820.6                 \\
1250           & 25            & 0        & 5                 & 0.3                     & 29480.0                 & 0                 & 609.8                   & 29680.0                 \\
1250           & 25            & 1        & 5                 & 3.6                     & 43983.0                 & 5                 & 144.3                   & 43688.2                 \\
1875           & 25            & 0        & 5                 & 0.4                     & 39740.0                 & 0                 & 618.7                   & 40260.0                 \\
1875           & 25            & 1        & 5                 & 5.6                     & 68745.8                 & 5                 & 239.3                   & 68088.0                
\end{tabular}
\caption{Comparison of runtime and solution quality between \textbf{CG} and \textbf{IP} on instances with express trains.}
\label{tab:express}
\end{table}

Table~\ref{tab:express} reports the results obtained.  The first three columns report the instance configuration.  The next three columns, $n_{\textbf{CG}}$, $\hat{t}_{\textbf{CG}}$, and $\hat{z}_{\textbf{CG}}$, report the number of instances solved, the average solution time, and the average objective function value for \textbf{CG}, respectively, over the five instances in that configuration.  The final three columns, $n_{\textbf{IP}}$, $\hat{t}_{\textbf{IP}}$, $\hat{z}_{\textbf{IP}}$, report the same statistics, but for \textbf{IP}.   All instances were solved by \textbf{CG} in reasonable time, i.e., under 6 seconds, whereupon \textbf{IP} is only able to solve those instances with $\alpha=1$ and even for these instances takes at least an order-of-magnitude more time, on average. Since \textbf{CG} is solving a more constrained version of the problem, its solutions can be worse than that of \textbf{IP}, but this improvement is marginal, at most $\sim 5\%$ worse.   Alternatively, when \textbf{IP} is unable to find the provably optimal solution, \textbf{CG} is able to find a better solution, up to 2$\%$ better than the solution obtained by \textbf{IP} in 10 minutes. Since the solutions time for \textbf{CG} are also often an order-of-magnitude faster than the solution times for \textbf{IP}, if an operator needs high-quality solutions quickly, the solution obtained by \textbf{CG} could suffice.

\subsubsection{Multiple destinations per CV trip} Another restrictive assumption inherent in the BDD-based approach is that the CVs visit one-and-only-one destination per trip.  \cite{iflmpicaps18} provided an analysis of how far from optimal the solutions obtained by single-destination-per-CV trip solutions are from those obtained by \texttt{IP} allowing passengers going to non-common destinations for small-scale problems.  We extend that analysis here and find even more encouraging results. 

We provide this comparison on the following instances. For this case, we used the same problem 
setting as in our conference paper~\cite{iflmpicaps18}.  For sake of brevity, we refer the interested 
reader to the section on Experiments in~\cite{iflmpicaps18} for a description on the routes for the 
commuter vehicles.  In this case, the number of destinations is 10.  
The number of passengers is chosen as $n \in \{500,750\}$.  
The train trips are retained as described in the beginning of the section.  
Again, we report results for $\alpha = 0,1$. 

For $\alpha = 0$, all instances were solved by \textbf{CG} in under one second, and none were solved by \textbf{IP} in over three hours.  \textbf{IP} only found a feasible solution to two instances in the three-hour time limit, with \textbf{CG} finding solutions at least as good as \textbf{IP} for both of those instances.  For $\alpha = 1$, \textbf{IP} ran out of memory, and \textbf{CG} solved all instances within 2.09 seconds.  Despite lacking optimality guarantees, \textbf{CG} can be used as a heuristic where \textbf{IP} fails to identify any feasible solution.  A significant drawback with the IP formulation 
is that the number of variables in the optimization problem scales with the number of 
possible route choices.  The number of possible route choices when passengers with different 
destinations share a CV grows exponentially as,
$\perm{K}{1}+\perm{K}{2}+\cdots+\perm{K}{\cvCap}$.
As a consequence, the loading of the \textbf{IP} model in memory consumes a significant 
amount of computational time ($\sim$1 hour) and solution of the linear relaxation at the 
root node also takes a comparable amount of time.  This emphasizes the need for developing alternative methods for this problem generalization.  Until such an algorithm is 
developed, \textbf{CG} offers a computationally inexpensive and 
scalable approach to obtaining high-quality solutions.

\section{Conclusion and future work}
\label{sec:conclusions}

In this paper we introduce an optimization algorithm for solving the problem of scheduling passengers on multiple legs of a last-mile transportation system. We study a version of the problem where in the last leg of the transportation system, passengers are transported via small-capacity commuter vehicles (CVs) to a limited set of destinations.  In particular, we study a variant of the problem where each CV trip carries passengers to a common destination, showing that this simplified version remains NP-hard.  The optimization framework relies on a decomposition of the problem into a collection of small-sized decision diagrams that can be mutually optimized over. This algorithm is shown to dramatically outperform existing techniques. Our experimental evaluation indicates that the algorithm developed can scale to problems of practical size and that considering conflicting objectives of minimizing passenger wait times with the total number of CV trips simultaneously does not significantly hinder the performance on either objective taken individually. 

The potential for expansion of this work is vast, as automated transportation networks becomes a reality.  The variant studied in this paper is a simplified version of real-world systems, where CVs can stop in multiple destinations.  The work 
of~\cite{iflmpicaps18} indicates that the gap between the optimal solutions obtained by limiting CVs to stop at only one destination per trip might not be far from the optimal solutions in the more general version of the problem. In fact, the solutions obtained by the decision diagram-based model can be used as a heuristic to more general variants, and the numerical evaluation in this paper suggests that the solutions can be found very quickly and are of high quality.  
Given the substantial savings in the long-run for any such improvement, adopting the methods developed in this paper to the more general case might be an interesting research direction.  Incorporating more real-world features like time-dependent travel times on the CVs, dynamic scheduling and response to traffic are additional dimensions that could potentially be added to the models developed in this paper. This work gives a critical and substantial step towards understanding how to solve challenging automated scheduling problems in the context of automated commuter systems. 


%
%
%
 \begin{APPENDICES}

\section{Proof of Theorem~\ref{thm:complexity}}
\label{sec:proofCompleixty}

\proof{Proof.}
We first show that the feasibility of ILMTP-SD is in NP. 
If we are given a solution consisting of the CV boarded by each passenger and the boarding time, 
then we can easily verify the feasibility of the solution. 
First, we check if, for each passenger, that there is a mass transit service that could bring the passenger to the terminal before the boarding time. In the worst case, this is proportional to the number of passengers times the number of mass transit trips. 
Second, for each CV we define a vector of tuples, each of which consists of the boarding time and destination of a passenger using that CV.
After sorting each of those vectors by the boarding times, 
we check with a linear 
pass on 
the vectors if 
(i) passengers boarding the CV at the same time have the same destination; and 
(ii) the time between consecutive trips is sufficient for the CV to return to the terminal. 

Next, we show that a decision version of the \emph{bin packing} problem can be reduced in polynomial time and space to the feasibility of the ILMTP-SD.  The bin packing problem can be stated as:  
Given a set of bins $B_1, B_2, \ldots$ with identical capacity $V$ and a list of $n$ items with sizes $a_1, \ldots, a_n$, does there exist a packing using at most $M$ bins?

The Karp reduction~\citep{Karp} is as follows. We define an ILMTP-SD instance with 
$n$ passengers and $M$ CVs, where passenger $p_i$ corresponds to item $i$ of the 
bin packing problem. Each of those passengers has a distinct destination $d_i$ and the round trip time 
is $\timeRoundTimeToBuilding{d_i} = a_i$. We assume that there is a single mass transit trip that can bring these passengers to $\terminal$ and that it arrives to $\terminal$ at time $t_0$, 
whereas the CVs must return to $\terminal$ by time $\maxTime = t_0 + V$. 
Finally, 
the origin of each passenger is irrelevant, the capacity of each CV is $\cvCap = 1$, we assume a time window $\timewindow = + \infty$, and the objective coefficient $\alpha = 1$. 

If there is a feasible solution to the ILMTP-SD problem above, 
then the bin packing problem has an affirmative answer. 
Namely, let $\mathcal{P}^i$ be the set of all passengers that board CV $i$ in the solution, 
which are on different trips since the CV capacity is 1. 
Assign the items corresponding to those passengers to bin $i$.  Since the first passenger in $\mathcal{P}^i$ to board CV $i$ left $\terminal$ after $t_0$ and the CV returned back to $\terminal$ before $\maxTime = V$, the sum of the duration of the trips for the passengers in $i$ must not exceed $V$. Therefore the associated items fit into bin $i$.
Hence, 
each CV corresponds to a bin and
all passengers boarding a given CV are assigned to that bin, using at most $M$ bins.

Conversely, if there is no feasible solution the ILMTP-SD problem, 
then the bin packing problem has a negative answer. 
If the bin packing problem has a solution using at most $M$ bins,  
then we can construct a solution for the corresponding ILMTP-SD instance 
through the same transformation. 

Therefore, the feasibility of the ILMTP-SD is as hard as the bin packing decision problem, 
which is known to be NP-complete~\citep{GJ}.
 \Halmos
\endproof

\section{Proof of Theorem~\ref{prop:optSolutionStructure}}\label{app:theorem2}

Before proceeding with the proof of Theorem~\ref{prop:optSolutionStructure}, we present a lemma.  

\begin{lem}
	\label{lem:1}
	Consider any feasible solution $\mathsf{g} = \set{g_1, \ldots, g_\gamma}$ to the ILMTP-SD , with 
	associated departure time from $\terminal$ of $t^{\mathsf{g}}_l$ , for $l = 1, \ldots, \gamma$.    
	For any two passengers $j', j''$ for which $\building(j') = \building(j'')$ and 
	$\mathsf{g}(j') \neq \mathsf{g}(j'')$, if the solution $\mathsf{g}' = \set{g'_1, \ldots, g'_\gamma}$ defined by
	\[g'_l = \left\{
	\begin{array}{ll}
	 g_l & : g_l \notin \set{\mathsf{g}(j'), \mathsf{g}(j'')}\\
	g_l \backslash \set{j'} \cup \set{j''} & : g'_l = \mathsf{g}(j') \\
	g_l \backslash \set{j''} \cup \set{j'} & : g'_l = \mathsf{g}(j''),
	\end{array}
	\right. \forall l = 1,\ldots,\gamma
	\]
	with $t^{\mathsf{g}'}_l = t^{\mathsf{g}}_l$, is a feasible solution, then the objective function value of both solutions is the same. 
\end{lem}

\noindent
\emph{Proof of Lemma~\ref{lem:1}. }
Switching the two passengers does not add or delete any CV trips, and hence, $(1-\alpha)$ times the number of CV trips remains unchanged.  We need only ensure that the total wait time at $\terminal$ remains unchanged. 

Let $l',l''$ be the indices of the groups that passengers $j',j''$ are assigned to in $\mathsf{g}$, respectively.  Furthermore, let $\tau', \tau''$ be the time that passengers $j',j''$ arrive to $\terminal$ on their mass transit trips, respectively.  The wait time at $\terminal$ for passenger $j'$ is changed from $t^{\mathsf{g}}_{l'} - \tau'$ to $t^{\mathsf{g}}_{l''} - \tau'$, a net change of $t^{\mathsf{g}}_{l''} - t^{\mathsf{g}}_{l'}$.  By the same argument, the net change in wait time for passenger $j''$ is $t^{\mathsf{g}}_{l'} - t^{\mathsf{g}}_{l''}$.  
The net changes in $j', j''$ cancel out.  Since wait time of other passengers are not affected, the result holds. 

Equipped with Lemma~\ref{lem:1}, we can now prove Theorem~\ref{prop:optSolutionStructure}.

\noindent
\emph{Proof of Theorem~\ref{prop:optSolutionStructure}. }
By way of contradiction, suppose there exists an instance for which there is no optimal solution 
satisfying the condition of the theorem for a $d^* \in \buildingSet$. Consider the optimal solution for which the smallest indexed passenger $j$ that violates this condition is maximized. 
Let $j^*$ be the smallest index in this solution for which there exists a $k$ with $\mathsf{g}(j^*) = \mathsf{g}(j^* + k)$ and $\mathsf{g}(j^* + 1 ) \neq \mathsf{g}(j^*)$.  Let $k^*$ be such an index $k$, and let $g_\ell = \mathsf{g}(j^*)$ and 
$g_{\ell'} = \mathsf{g}(j^*+1)$.
We first show that  
$j \geq j^*+1$ for all $\{j \,|\, j \in {g}_{\ell'}\}$.  Suppose not;  let $\hat{j} = \arg\max \{ j \, | \, 
j \in {g}_{\ell'}, j < j^*+1  \}$ (which will be non-empty by assumption).  
Then passenger $j'$, $j'+k'$ with $j' = \hat{j}, k' = j^*+1-\hat{j}$ satisfy 
$j', j'+k' \in \mathsf{g}_{\ell'}$ and 
$j'+1 \notin \mathsf{g}_{\ell'}$ are a set of indices violating the claim of the theorem with $j' < j^*$, 
contradicting the minimality of $j^*$. 
Hence, $(j^*+1)$ is the minimum index among all $j \in g_{\ell'}$.  

In the remainder of the proof, we construct another solution in which the index $j^*$ does not violate 
the claims of the theorem.  This 
contradicts the maximality of $j^*$ among all optimal solutions, thereby establishing the result.  

Conditioning on the relative values of the departure times of the CVs for ${g}_\ell$ and ${g}_{\ell'}$, first consider the case where $t^{\mathsf{g}}_\ell \leq t^{\mathsf{g}}_{\ell'}$.  We claim that exchanging the group assignment of $j^*+1$ and $j^*+k^*$ and holding all else equal results in another feasible solution.   For all passengers $j$, let $\rel(j), \ded(j)$ be $\passengerTimeRequest{j} - \timewindow,\passengerTimeRequest{j} + \timewindow$, respectively.  We need only show that (a) $ \rel(j^*+1) \leq t^{\mathsf{g}}_{\ell}+\timeTerminalToBuilding{d^*} \leq \ded(j^*+1) $  and (b) $ \rel(j^*+k^*) \leq t^{\mathsf{g}}_{\ell'}+\timeTerminalToBuilding{d^*} \leq \ded(j^*+k^*) $.  (a) follows because any passenger can be placed in a group with fewer than $\cvCap$ passengers  without breaking feasibility if the passenger's request time is before or after at least one other passenger.  The first inequality in (b) follows because  $ \rel(j^*+k^*) \leq t^{\mathsf{g}}_{\ell}+\timeTerminalToBuilding{d^*}$, by the feasibility of the original solution, and  $t^{\mathsf{g}}_{\ell}+\timeTerminalToBuilding{d^*} \leq t^{\mathsf{g}}_{\ell'}+\timeTerminalToBuilding{d^*}$, by assumption.  The second inequality in (b) holds because $\ded(j^*+1) \leq t^{\mathsf{g}}_{\ell'}+\timeTerminalToBuilding{d^*}$, by 
the feasibility of the original solution and $\ded(j^*+1) \leq \ded(j^*+k^*)$, by the ordering 
of passengers.   Furthermore, by Lemma~\ref{lem:1} the objective function remains unchanged by this 
exchange, and is therefore optimal.  If the resulting solution satisfies the claim of this theorem, then 
the claim holds.  
If not, then the claim of this theorem is violated for another $j > j^*$;  contradicting the 
maximality of $j^*$. 

We now consider the alternative case where 
$t^{\mathsf{g}}_{\ell} > t^{\mathsf{g}}_{\ell'}$. The exchange from the previous case may not work because putting passenger $j^* + k^*$ into group ${g}_{\ell'}$ may not be feasible.  Additionally, if there exist $k' > 0$ passengers in $g_\ell$ with indices lower than $j^*$ then these passengers 
must be $j^*-k',\ldots,j^*-1$.  If not, it would contradict the assumption of $j^*$  
as the smallest index in the optimal solution violating the claim of this theorem. Define $k'$ so that $j^* - k'$ is the minimum indexed passenger in group ${g}_\ell$, which is 0 if $j^*$ is the minimum indexed passenger.

Consider the following two-step exchange\textemdash for $i = 0, \ldots, k'$, move each passenger $j^*-i$ 
from ${g}_{\ell}$ into ${g}_{\ell'}$. Then, move the $k'+1$ passengers with the highest indices in the resulting ${g}_{\ell'}$ into ${g}_\ell$.  The resulting groups have the same cardinality as they originally had, and so by Lemma~\ref{lem:1} the objective values remain the same.

We now show that the resulting solution is valid and then show that the choice of optimal solution contradicts the maximality among optimal solutions assumption on the selection of $j^*$, which concludes the proof.  Any passenger $j^*-i \in {g}_\ell$ for $i = 0,...,k'$ can be moved to ${g}_{\ell'}$ without violating $(j^*-i)$'s arrival time window because $\rel(j^*-i) \leq \rel(j^*+1) \leq t^{\mathsf{g}}_{\ell'} +\timeTerminalToBuilding{d^*} < t^{\mathsf{g}}_{\ell}+\timeTerminalToBuilding{d^*}$ and  $\ded(j^*-i) \geq t^{\mathsf{g}}_{\ell}+\timeTerminalToBuilding{d^*} > t^{\mathsf{g}}_{\ell'}+\timeTerminalToBuilding{d^*}$.  Additionally, any passenger $j \in {g}_{\ell'} $ can be moved to ${g}_\ell$ without violating the arrival time windows because $\rel(j) \leq t^{\mathsf{g}}_{\ell'}+\timeTerminalToBuilding{d^*} < t^{\mathsf{g}}_{\ell}+\timeTerminalToBuilding{d^*}$ and 
$\ded(j) \geq \ded(j^*+1) \geq \ded(j^*) \geq t^{\mathsf{g}}_{\ell}+\timeTerminalToBuilding{d^*}$ where first inequality follows from the result that $(j^*+1)$ is the minimum index for all $j \in g_{\ell'}$.  
Hence, the resulting solution is also optimal.  Finally, if the 
resulting solution violates the claim of this theorem, then the smallest index must be larger than $j^*$.  
This again contradicts the maximality of $j^*$ among all optimal solutions, as assumed. 

\section{Proof of Theorem~\ref{theorem:DDproperty}}\label{app:DDproperty}

\proof{Proof.}

Constructivelly from Algorithm~\ref{alg:dd}, the number of nodes in the DD is at most $(\nPassengersPerBuilding{d} + 1) \cdot \cvCap = O(\nPassengersPerBuilding{d} \cdot \cvCap)$.  Each node has at most $\timewindow\cdot2 + 1$ arcs, so that the maximum number of arcs is $(\nPassengersPerBuilding{d} + 1) \cdot \cvCap \cdot (\timewindow*2 + 1) = O(\nPassengersPerBuilding{d} \cdot \cvCap \cdot \timewindow)$, as desired.  Since it takes constant time to create each arc, the time bound follows.  

For property~\textnormal{\textbf{(DD-1)}}, fix an arbitrary destination $d$, and consider any path $p \in \pathSet^d$ and a passenger $j^d_i$.  Consider the first one-arc in $p$ below layer $i$ (including that layer).  Note that since the last layer only contains one-arcs, such an arc exists.  Suppose this one-arc $a'$ is directed out of node $\DDnode$ in layer $L^d_{i'}$, with $i' \geq i$.  By construction, the state of $\DDnode$ must be greater than or equal to $i' - i$.  $\group(a)$ therefore contains $j^d_i$, and so $j^d_i$ is in at least one group in $\groups{p}$.  Furthermore, any one-arc below layer $L^d_{i'}$ that connects $\arcterminal{a'}$ to $\terminalnode^b$ can only select passengers $j^d_{i''}$, with $i'' \geq i'$, and so $j^d_i$ appears in exactly one $\groups{p}$.

Property~\textnormal{\textbf{(DD-2)}} is satisfied by construction\textemdash each one-arc $a'$ directed out of node $\DDnode$ satisfies that every passenger in $g(a)$ will arrive to $d$ within the desired 
range of arrival times.  

What remains to be shown is that property~\textnormal{\textbf{(DD-3)}} is satisfied.   Fix an arbitrary $\mathsf{g}' \in \mathcal{G}^d$.  Let $g_1, \ldots, g_f$ be the $f$ groups in $\mathsf{g}'$, ordered by the index of the passengers in the groups. Let $j^d_{i_h}$ be the highest indexed passenger in each group $h$, $h=1, \ldots, f$.  We proceed by induction on $h$ to show that there is a sequence of zero-arcs from the node on layer $\DDlayer^d_{i_{h}}$ to layer $\DDlayer^d_{i_{h+1}-1}$ to a node $\DDnode$ with state $i_{h+1}-1 - i_h$, and then a one-arc from $\DDnode$ to the node on layer $\DDlayer^d_{i_{h+1}}$ with state 0 for every time that the passengers can mutually depart from $\terminal$. 

We first establish the base case.  Starting from $\rootnode^d$, follow $i_1-1$ zero-arcs until layer $\DDlayer^d_{i_1-1}$.  Because $g_1$ is in $\mathsf{g}'$, each $t \in \set{t^e(j^d_{i_1}), \ldots, t^l(j^d_1)}$ is a feasible departure time from $\terminal$ for $g_1$, and so a one-arc directed out of the node on this layer with state $i_1 - 1$ will be added to $\DD^d$ with the arc-start-time $t$. 

By induction on $h$, consider the node $\DDnode$ on layer $\DDlayer^d_{i_h}$ with state 0.  Follow $i_{h+1}-1-i_h$ zero-arcs.  This will end at a node $\DDnode'$ on layer $\DDlayer^d_{i_{h+1}-1}$ with state $i_{h+1}-1-i_h$.  Since $g_h \in \mathsf{g}'$, each time $t \in \set{t^e\left(j^d_{i_{h+1}}\right), \ldots, t^l\left(j^d_{i_{h}}+1\right)}$ is a feasible departure time from $\terminal$ for all passengers in $g_{j}$, and so a one-arc from $\DDnode'$ to the node on layer $\DDlayer^d_{i_{h+1}}$ with state 0 will be added to $\DD^d$ with this arc-start-time.
 


 \end{APPENDICES}


\bibliographystyle{informs2014} 
\bibliography{iflmp} 


\end{document}